\documentclass[letterpaper,11pt]{article}%
\usepackage{amsmath,amssymb,amsthm,latexsym}
\usepackage{hyperref}
\usepackage{epsfig,subfigure}
\usepackage{fullpage}
\usepackage{amsmath}
\usepackage{amsfonts}
\usepackage{amssymb,bbold}
\usepackage{graphicx,color}%
\usepackage[margin=1.15in]{geometry}
\usepackage{enumerate}
\usepackage[cmyk]{xcolor}
\usepackage{mathtools}
\usepackage{pdfsync}
\usepackage{dsfont}
\usepackage{color}
\usepackage{bm}
\usepackage{cite}
\usepackage{bbm}
\usepackage{multirow}
\usepackage{tabularx,ragged2e,booktabs,caption}

\providecommand{\U}[1]{\protect\rule{.1in}{.1in}}

\theoremstyle{plain}
\newtheorem{thm}{Theorem}[section]
\newtheorem{prop}[thm]{Proposition}
\newtheorem{lem}[thm]{Lemma}

\theoremstyle{remark}
\newtheorem{remark}[thm]{Remark}

\theoremstyle{definition}
\newtheorem{defn}{Definition}[section]

\newcommand{\e}{\varepsilon}
\newcommand{\disp}{\displaystyle}

\newcommand{\ba}{\begin{array}}
\newcommand{\ea}{\end{array}}

\newcommand{\bthm}{\begin{thm}}
\newcommand{\ethm}{\end{thm}}
\newcommand{\bprop}{\begin{prop}}
\newcommand{\eprop}{\end{prop}}
\newcommand{\blemma}{\begin{lem}}
\newcommand{\elemma}{\end{lem}}

\newcommand{\beqn}{\begin{equation}}
\newcommand{\eeqn}{\end{equation}}
\newcommand{\beqns}{\begin{equation*}}
\newcommand{\eeqns}{\end{equation*}}

\newcommand{\supp}{\operatorname{supp}}

\newcommand{\pr}{\prime}
\newcommand{\pt}{\partial}
\newcommand{\arrow}{\rightarrow}

\newcommand{\warrow}{\rightharpoonup}

\newcommand{\Prob}{\mathcal{P}}

\renewcommand{\leq}{\leqslant}
\renewcommand{\geq}{\geqslant}

\newcommand{\Om}{\Omega}

\definecolor{mygreen}{rgb}{0.1,0.75,0.2}

\newcommand{\rhu}{\rightharpoonup}

\newcommand{\E}{\mathsf{E}}
\newcommand{\Se}{\mathsf{S}}

\newcommand{\tbf}{\mathbf{t}}
\newcommand{\ibf}{\mathbf{i}}
\newcommand{\xibf}{\bm{\xi}}

\newcommand{\Rd}{{\mathord{\mathbb R}^d}}
\newcommand{\loc}{{\rm loc}}

\def\P{{\mathcal P}}

\newcommand{\scie}{{\textrm{e}}}

\makeatletter
\newcommand{\leqnomode}{\tagsleft@true\let\veqno\@@leqno}
\newcommand{\reqnomode}{\tagsleft@false\let\veqno\@@eqno}
\makeatother

\newcommand{\R}{{\mathbb R}}

\newcommand{\One}{\mathbb{1}}

\newcommand{\munubar}{\bar{\mu}_\nu}

\newcommand{\Pa}{{\mathcal{P}_2^a}}

\newcommand*\samethanks[1][\value{footnote}]{\footnotemark[#1]}

\begin{document}

\title{Swarming in domains with boundaries: approximation and regularization by nonlinear diffusion}

\author{R.C. Fetecau \thanks{Department of Mathematics, Simon Fraser University, 8888 University Dr., Burnaby, BC V5A 1S6, Canada}
\and M. Kovacic \samethanks
\and I. Topaloglu \thanks{Department of Mathematics and Applied Mathematics, Virginia Commonwealth University, 1015 Floyd Ave, Richmond, VA 23219, United States}
}

\maketitle

\begin{abstract}
We consider an aggregation model with nonlinear diffusion in domains with boundaries and investigate the zero diffusion limit of its solutions. We establish the convergence of weak solutions for fixed times, as well as the convergence of energy minimizers in this limit. Numerical simulations that support the analytical results are presented. A second key scope of the numerical studies is to demonstrate that adding small nonlinear diffusion rectifies a flaw of the plain aggregation model in domains with boundaries, which is to evolve into unstable equilibria (non-minimizers of the energy).
\end{abstract}


\section{Introduction} \label{sec:intro}

In this paper we investigate the zero diffusion limit ($\nu\to 0$) of weak solutions to the following aggregation equation:
	\beqn \label{eqn:eps-eqn}
		\left \{
			\begin{aligned}
				\pt_t \mu_\nu + \nabla \cdot (\mu_\nu \mathbf{v}_\nu ) &= 0 \qquad \text{in } \Om \times [0,T], \\
				  \text{with} \quad \mathbf{v}_\nu &= -\frac{\nu^\alpha \, m}{m-1} \nabla \rho_\nu^{m-1} - \nabla K * \mu_\nu - \nabla V, \\
				 \mathbf{v}_\nu \cdot \mathbf{n} &= 0 \qquad \text{on } \pt\Om \times [0,T], \\
				  \mu_\nu (0) &= \mu^0 \qquad \text{on } \Om,
			\end{aligned} 
		\right.
	\eeqn
where $\mu_\nu$ is an absolutely continuous probability measure with density $\rho_\nu$, $K$ is an interaction potential, and $V$ is an external potential. The equation is set in a closed domain $\Om\subset\Rd$ with smooth boundary, and $\mathbf{n}$ denotes the outward unit normal to $\pt\Om$. Also,  $\nu>0$ is the diffusion coefficient (with exponent $\alpha>0$) and $m>1$.

The aim of the paper is to study how solutions of \eqref{eqn:eps-eqn} approximate solutions of the first-order aggregation model given by: 
	\beqn \label{eqn:agg-eqn2}
		\left \{
			\begin{aligned}
				\pt_t \mu + \nabla \cdot (\mu \mathbf{v} ) &= 0 \qquad \text{in } \Om \times [0,T], \\
				 \text{with} \quad \mathbf{v} &= P_x(- \nabla K * \mu - \nabla V), \\
			 	\mu (0) &= \mu^0 \qquad \text{on } \Om,
			\end{aligned} 
		\right.
	\eeqn
where $P_x:\R^d \to \R^d$ is a projection operator defined by
		\beqn \label{eqn:proj-op}
			P_x \xi \,=\, \begin{cases*} 
											\xi &   if $x\not\in\pt\Om$ or if $x\in\pt\Om$ and $\xi \cdot \mathbf{n} \leq 0$, \\
										    \Pi_{\pt\Om}\xi &  otherwise.
									    \end{cases*}
	\eeqn	
Here $\Pi_{\pt\Om}$ is the projection onto the tangent plane of the boundary.  Note that in model \eqref{eqn:agg-eqn2}, $\mathbf{v} = - \nabla K * \mu - \nabla V$ everywhere in $\Om$, except at points $x \in \partial \Om$ where $- \nabla K * \mu - \nabla V$ points {\em outward} the domain, in which case its projection on the tangent plane to the boundary is considered instead.  This is a ``slip, no-flux" boundary condition by which particles/individuals that meet the boundary, do not exit the domain, but move freely along it \cite{WuSlepcev2015, CarrilloSlepcevWu2016}.  Throughout this work we will be referring to model \eqref{eqn:agg-eqn2} as the {\em plain} aggregation model.


Model \eqref{eqn:agg-eqn2}, typically set in free space with no boundaries ($\Omega = \R^d$), has been a topic of intense interest in recent years. The model appeared in various contexts related to swarming and social aggregations, such as biological swarms and pattern formation \cite{M&K, KoSuUmBe2011, LeToBe2009}, granular media~\cite{Toscani2000, CaMcVi2006}, self-assembly of nanoparticles~\cite{HoPu2005}, Ginzburg-Landau vortices~\cite{DuZhang03},  robotics and space missions \cite{JiEgerstedt2007} and opinion dynamics \cite{MotschTadmor2014}. There is an extensive literature on the mathematical properties of the first-order model in free space, which includes studies on the well-posedness of solutions \cite{BodnarVelasquez2, BertozziLaurent, BertozziCarilloLaurent}, the long-time  behaviour of solutions \cite{FeHuKo11, FeHu13, LeToBe2009}, and blow-up (in finite or infinite time) by mass concentration \cite{FeRa10, BertozziCarilloLaurent, HuBe2010}. 	

Practical setups of model 	\eqref{eqn:agg-eqn2} involve domains with nontrivial boundaries (e.g., the model for locust dynamics in  \cite{ToDoKeBe2012}, or applications in environments with impenetrable walls or obstacles). Despite its high relevance to applications however, model \eqref{eqn:agg-eqn2} in domains with boundaries has been considered in only a few works \cite{BeTo2011, WuSlepcev2015, CarrilloSlepcevWu2016}. Most relevant to the present paper is the study in \cite{FeKo2017} where authors identified a degeneracy of model \eqref{eqn:agg-eqn2}, namely that its solutions tend to evolve into unstable equilibria. This is an intriguing flaw of the plain aggregation model, in particular since it has a gradient flow formulation. Investigating the nonlinear diffusion model \eqref{eqn:eps-eqn} (in the zero diffusion limit) as a means to rectify this degeneracy of the plain aggregation model, is the main motivation for the present research.

Nonlinear diffusion models of class \eqref{eqn:eps-eqn} have also received a great deal of interest lately. In the context of modelling biological aggregations, model \eqref{eqn:eps-eqn}  (with $m=3$) was discussed in \cite{TBL}; of particular interest for such applications is that the model can exhibit equilibria with compact support and sharp edges that correspond to localized clumps of organisms. There is also extensive work on mathematical studies of these models (mostly set in free space), including well-posedness results\cite{AGS2005, BertozziSlepcev2010, BeRoBe2011}, investigations into properties of the steady states \cite{BurgerDiFrancescoFranek, HuFe2013, BuFeHu14, Kaib17, CaHiVoYa16}, and studies of long-time behaviour of the solutions\cite{Burger:DiFrancesco, CaMcVi2006, ChKiYa2013}. 
In addition, linear diffusion models have an extensive literature on their own. Relevant to the ideas of this paper, \cite{Zh2017} studies the zero diffusion limit of the linear diffusion model, and  \cite{EvKo16} showcases how linear diffusion can remove the degeneracy of multiple equilibria and lead to a unique steady state.

The major goal of the present work is to show that model \eqref{eqn:eps-eqn} approximates, as well as {\em regularizes}, the plain aggregation model  \eqref{eqn:agg-eqn2}. Regarding the former, we establish several results. First, we show that at each {\em fixed} $t>0$, weak solutions of the diffusive model \eqref{eqn:eps-eqn} converge to solutions of \eqref{eqn:agg-eqn2} in the zero diffusion limit. Second, we study the convergence of minimizers of the energy associated to model \eqref{eqn:eps-eqn} to minimizers of the energy of the plain aggregation model (more details in the methodology and results presented below). On the other hand, concerning the regularization thesis, we provide strong numerical evidence that the diffusive model rectifies the aforementioned flaw of model \eqref{eqn:agg-eqn2}, in the sense that solutions of \eqref{eqn:eps-eqn} bypass the unstable equilibria of model \eqref{eqn:agg-eqn2}. Note however that, in consistency with the approximation result at fixed times, the solutions of the diffusive model can pass arbitrarily close to this equilibrium, as diffusion becomes arbitrarily small. 

\paragraph{Methodology and main results.} The mathematical analysis of models \eqref{eqn:eps-eqn} and \eqref{eqn:agg-eqn2}, as investigated in numerous papers, is based on the formulation of these models in terms of gradient flows on spaces of probability measures equipped with the Wasserstein metric, using the theory developed in \cite{AGS2005}. The diffusive model is discussed in \cite[Chapter 11]{AGS2005}, while model  \eqref{eqn:agg-eqn2} was studied more recently in \cite{WuSlepcev2015, CarrilloSlepcevWu2016}. Specifically, models \eqref{eqn:eps-eqn} and \eqref{eqn:agg-eqn2} have a gradient flow structure and their equilibria are stationary points of the following energy functionals, respectively: 	\beqn \label{eqn:reg_energy}
		\E_\nu(\mu) \,=\, \begin{cases*} 
											\frac{\nu^\alpha}{m-1} \int_{\Om} \rho^m(x)\,dx + \frac{1}{2}\int_{\Om}\!\int_{\Om} K(x-y)\rho(x)\rho(y)\,dxdy & \  if $d\mu=\rho\,dx$ \\
											\qquad\qquad\qquad\qquad\qquad\qquad\qquad + \int_{\Om} V(x)\rho(x)\,dx & \   with $\rho\in L^1\cap L^m$ \\
											\\
										    + \infty & \quad otherwise,
									    \end{cases*}
	\eeqn
and 
	\beqn\label{eqn:energy}
		\E(\mu) \, = \,  \frac{1}{2} \int_{\Om}\!\int_{\Om} K(x-y)\,d\mu(x) d\mu(y) + \int_{\Om} V(x)\,d\mu(x).
	\eeqn	
The two energies are defined over $\P_2(\Om)$, the space of probability measures with finite second moments.

The following general conditions on the interaction potential $K$, external potential $V$, and domain $\Om$ will be referenced in the analytical results presented below. \leqnomode
		\begin{itemize}	
				\addtolength{\itemsep}{6pt}
						\item[(H1)] $K$ is symmetric, i.e., $K(x)=K(-x)$ for all $x\in\Rd$.
						\item[(H2)] Either
												\beqn  \label{eqn:KassmA} \tag{H2A}
												\begin{gathered}
													K\in C^2(\Om), \ K \text{ is bounded from below}, 													
													\text{and } K \text{ is } \lambda\text{-convex for some } \lambda\in\R,
												\end{gathered}
												\eeqn
											or
												\beqn \label{eqn:KassmB} \tag{H2B}
												\begin{gathered} 
													K \in L^1_\loc(\Om) \cap C^1(\Om \setminus B_R(0)), \
														K \to +\infty \text{ as } |x|\to \infty, \\
													\text{there exists a function } K^a\in C(\Om) \text{ so that } K^r:=K-K^a \text{ is superharmonic}, \\
													\text{ and for } |x| \text{ large, } K^r \text{ and } K^a \text{ have at most quadratic growth}.
												\end{gathered}
												\eeqn
						\item[(H3)] Either
												\beqn \label{eqn:VassmA} \tag{H3A}
														V\in C^1(\Om), \ |\nabla V(x)| \leq C(1+|x|) \text{ for all } x\in\Rd, \text{ and } V \text{ is } \lambda\text{-convex},
												\eeqn
										or
												\beqn \label{eqn:VassmB} \tag{H3B}
													\begin{gathered} V\in C(\Om), V \text{ has at most quadratic growth, it is bounded from below, and} \\
														V \text{ is either strictly increasing or translation invariant in every component}.
													\end{gathered}
												\eeqn												
						\item[(H4)] $\Om\subset\Rd$ is bounded, convex, and $\pt\Om \in C^1$.
		\end{itemize}
\reqnomode

The conditions that appear with the extra label A or B are used separately, for the convergence of weak solutions and for the variational results, respectively. Condition (H4) is assumed only for the dynamic convergence results.

Our first main result regarding the convergence of weak solutions is the following theorem.
\bthm \label{thm:conv}
	Assume $K$, $V$, and $\Om$ satisfy the hypotheses \emph{(H1), (H2A), (H3A)}, and \emph{(H4)}. Suppose $\mu_\nu(t)$ is a weak solution of \eqref{eqn:eps-eqn}, and $\mu(t)$ is a weak solution of \eqref{eqn:agg-eqn2} for all $t\in[0,T]$. Then there exist constants $C,\, \tilde{C} >0$ such that
		\[
			d_W^2(\mu_\nu(t),\mu(t)) \leq \tilde{C} \nu^{\beta} \,t\,e^{C t}
		\]
for all $t\in [0,T]$, where $\beta = \min \big\{\alpha-dm/(d+2),\, 1/(d+2)\big\}$, $\alpha>dm/(d+2)$, and $d_W$ denotes the 2-Wasserstein metric.
\ethm

Our proof follows the proof of \cite[Theorem 4.3]{Zh2017} where Zhang obtained a similar result in the case of linear diffusion and time-dependent domains. The idea here is to use the fact that, under stronger assumptions on $K$ and $V$, such as $\lambda$-convexity, and under boundedness and convexity assumptions on the domain $\Om$, weak solutions of the equation \eqref{eqn:eps-eqn} are gradient flow solutions of the energy functional \eqref{eqn:reg_energy}.

Our second result establishes the convergence of minimizers of diffuse-level energies $\E_\nu$ as the diffusitivity constant $\nu$ approaches zero. Compared to Theorem \ref{thm:conv}, the result applies to more general potentials $K$, external potentials $V$, and does not assume boundedness or convexity of the domain (assumption (H4)). However, we require a slightly stronger assumption on the exponent $\alpha$ than in Theorem \ref{thm:conv}. 

\bthm \label{thm:conv_min}
Assume $K$ and $V$ satisfy the hypotheses \emph{(H1), (H2B)}, and \emph{(H3B)}. Also assume $\alpha > dm$. Let $\{\mu_\nu\}_{\nu>0} \subset \P_2(\Rd)$ be a sequence such that each $\mu_\nu$ is a minimizer of the energy $\E_\nu$. Then there exists $\mu \in \P_2(\Rd)$ such that, up to a subsequence, $\mu_\nu \to \mu$ in $\P_2(\Rd)$ as $\nu\to 0$, and $\mu$ is a minimizer of $\E$. 
\ethm

In order to prove this theorem we utilize a standard technique in the calculus of variations, referred to as $\Gamma$-convergence. Theorem \ref{thm:conv_min} follows from the $\Gamma$-convergence of the sequence of energy functionals $\{\E_\nu\}_{\nu>0}$ combined with compactness of a sequence of minimizers of $\E_\nu$. 

We point out that Theorem \ref{thm:conv_min} applies to a wide range of interaction potentials $K$, which might have locally integrable singularities and which are not necessarily $\lambda$-convex. Indeed,  assumptions (H1) and (H2B)  are in particular satisfied by interaction potentials in attractive-repulsive power-law form: $K(x) = |x|^q/q - |x|^p/p \;$ for $2-d \leq p< q \leq 2$; this includes for instance singular Newtonian repulsion ($p = 2-d$) in $d\geq 3$. Also,  in 
dimensions $d \geq 2$, the assumptions are satisfied by interaction potentials with power-law attraction and logarithmic repulsion: $K(x) = |x|^q/q - \log(|x|)$ for $0< q \leq 2$; this includes Newtonian repulsion in $d=2$.

The summary of the paper is as follows. In Section \ref{sec:background} we provide preliminaries and background. Section \ref{sec:analysis} contains the proofs of Theorems \ref{thm:conv} and \ref{thm:conv_min}. In Section \ref{sect:numerics} we present numerical support for Theorems \ref{thm:conv} and \ref{thm:conv_min}, as well as for the regularization thesis, in one dimension.

\section{Preliminaries and background} \label{sec:background}

\subsection{The 2-Wasserstein metric}
We consider the space
		\[
			\P_2(\Om) := \left\{ \mu\in\P(\Om) \colon \int_{\Om} |x|^2\,d\mu(x) < +\infty \right\}
		\] 
of probability measures on $\Om$ with finite second moments, endowed with the 2-Wasserstein metric.  We recall this space briefly below, along with some of its basic properties.  For further background, we refer the reader to the books by Ambrosio, Gigli and Savar\'e \cite{AGS2005} and Villani \cite{Villani}.

The 2-Wasserstein distance between $\mu,\sigma \in \P_2 (\Om)$ is given by
		\beqn\label{eqn:Wass_dist}
			d_W(\mu,\sigma):=\left(\min\left\{\int_{\Om}\!\int_{\Om}|x-y|^2\,d\gamma(x,y)\colon \gamma\in\Gamma(\mu,\sigma)\right\}\right)^{1/2},
		\eeqn
where $\Gamma(\mu,\sigma)$ is the set of transport plans between $\mu$ and $\sigma$,	
		\[
			\Gamma(\mu,\sigma):=\left\{\gamma\in\P(\Om\times\Om)\colon (\pi_1)_\#\gamma=\mu\quad\text{and}\quad(\pi_2)_\#\gamma=\sigma\right\}.
		\]
Here $\pi_1$, $\pi_2$ denote the projections $\pi_1(x,y)=x$ and $\pi_2(x,y)=y$. For $i=1,2$, $(\pi_i)_{\#}\gamma$ denotes the pushforward of $\gamma$ defined by $(\pi_i)_{\#}\gamma (U):=\gamma(\pi_i^{-1}(U))$ for any measurable set $U\subset\Om$.

The minimization problem \eqref{eqn:Wass_dist} admits a solution, i.e., there exists an optimal transport plan $\gamma_0\in \Gamma(\mu,\sigma)$ so that
		\[
			d_W^2(\mu,\sigma)=\int_{\Om}\!\int_{\Om} |x-y|^2\,d\gamma_0(x,y).
		\]
If $\sigma$ is absolutely continuous with respect to the Lebesgue measure, then there is an optimal transport map $\tbf_\sigma^\mu : \Om \to \Om$ that transports $\sigma$ onto $\mu$ (i.e., ${\tbf_\sigma^\mu}_{\#}\sigma = \mu$) such that (see \cite{McCannExistence} for details):
		\[
			d_W^2(\mu,\sigma) = \int_{\Om} \big| \tbf_\sigma^\mu(x) - x \big|^2\,d\sigma(x).
		\]

Moreover, $(\P_2(\Om),d_W)$ is a complete and separable metric space, and convergence in $(\P_2(\Om),d_W)$ can be characterized as follows:
\smallskip
	\begin{center}
			\begin{tabular}{lcl}
				$d_W(\mu_n,\mu)\arrow 0$ & $\Longleftrightarrow$ & $\mu_n\arrow\mu$ weak-$*$ in $\P(\Om)$ and $\int_{\Om}|x|^2\,d\mu_n(x)\arrow\int_{\Om}|x|^2\,d\mu(x)$,
				\\
				\\
					& $\Longleftrightarrow$  & $\int_{\Om} f(x)\,d\mu_n(x) \to \int_{\Om} f(x) \,d\mu(x)$,\\[0.25cm]
					&   & $\forall f\in C(\Om)$ such that $|f(x)| \leq C(1+ |x-x_0|^2)$.
			\end{tabular}
	\end{center}
\medskip
We will refer to functions satisfying $|f(x)| \leq C(1+ |x-x_0|^2)$, for some $C>0$ and $x_0 \in \Om$, as functions with \emph{at most quadratic growth}.

Denote by $\Pa(\Om)$ the space of probability measures in $\P_2(\Om)$ that are absolutely continuous with respect to the Lebesgue measure. Given $\mu_1,\, \mu_2 \in \P_2(\Om)$, and $\sigma\in\Pa(\Om)$, we denote by $d_\sigma$ the \emph{pseudo-Wasserstein} distance with base $\sigma$ defined as
	\[
		d_\sigma^2(\mu_1,\mu_2) := \int_{\Om} \big| \tbf_\sigma^{\mu_1}(x) - \tbf_\sigma^{\mu_2}(x) \big|^2 \, d\sigma(x).
	\] 
By \cite[Proposition 1.15]{Craig16}, $d_\sigma$ defines a metric on $\P_\sigma(\Om) = \{ \mu\in \P(\Om) \colon d_W(\mu,\sigma) < +\infty \}$, and $d_W (\cdot,\cdot) \leq d_\sigma(\cdot,\cdot)$ for any $\sigma\in\Pa(\Om)$.

As shown by Ambrosio, Gigli and Savar\'{e}, the map $\mu \mapsto d_W^2(\mu,\sigma)$ is not convex on geodesics in $\P_2(\Om)$ but it is convex along generalized geodesics \cite[Lemma 9.2.1]{AGS2005}. Generalized geodesics whose base are in $\Pa(\Om)$ are constant speed geodesics with respect to the pseudo-Wasserstein metric \cite[Proposition 1.15]{Craig16}. Due to its connection with convexity properties of $\E_\nu$ along generalized geodesics, this metric also appears in HWI-type estimates controlling energies $\E_\nu$ and $\E$ along geodesics. In the proof of our first main result, Theorem \ref{thm:conv}, we control the distance between $\mu_\nu$ and $\mu$ which requires a decay control on $\mu$. As Zhang points out in \cite[Appendix A]{Zh2017}, a simple convolution does not provide sufficient decay rates in the pseudo-Wasserstein distance. In order to overcome this problem we utilize the following modification lemma.


\blemma(see Lemma 4.9 in \cite{Zh2017}) \label{lem:modify}
	For any $s>0$ sufficiently small, and for any $\mu\in \P_2(\Om)$, $\sigma\in \Pa(\Om)$ where $\Om\subset \Rd$ satisfies the hypothesis (H4), there exists $\mu_s\in\Pa(\Om)$ with density $\rho_s$ such that
		\[
			d_\sigma (\mu,\mu_s) \leq C_{\Om} \, s \qquad \text{and} \qquad \|\rho_s\|_{L^\infty(\Om)} \leq s^{-d}.
		\]
		The constant $C_\Om$ only depends on the domain $\Om$.
\elemma

\medskip

\subsection{Definitions of weak solutions} Next we define the notion of \emph{weak solutions} for the equations \eqref{eqn:eps-eqn} and \eqref{eqn:agg-eqn2}. To this end, recall that a curve $\mu(t): (a,b) \to \P_2(\Om)$ is \emph{locally 2-absolutely continuous} if there exists $m \in L^2_{\loc}(a,b)$ so that 
		\[
			d_W(\mu(t),\mu(s)) \leq \int_s^t m(r)\,dr \text{ for all } a <s \leq t < b .
		\]

\begin{defn}
For $\mu^0\in\Pa(\Om)$, a locally 2-absolutely continuous curve $\mu_\nu:[0,T]\to \Pa(\Om)$ is a weak solution of \eqref{eqn:eps-eqn} if for all $\varphi \in C_c^\infty(\Om \times (0,T))$,
 \[
				\int_0^T \! \int_{\Om} \big(\pt_t \varphi + \nu^\alpha \rho_\nu^{m-1}\Delta \varphi - \nabla V \cdot \nabla \varphi - (\nabla K * \mu_\nu)\cdot \nabla \varphi \big) \,d\mu_\nu \, dt = 0, \qquad\quad \mu_\nu(0)=\mu^0,
			\] 
and for all $t\in[0,T]$, $\supp \mu_\nu(t) \subset \Om$ where $\rho_\nu$ denotes the density of $\mu_\nu$ (see \cite[Section 11.2]{AGS2005}).
\end{defn}

\begin{defn}
For $\mu^0\in \P_2(\Om)$, a locally 2-absolutely continuous curve $\mu:[0,T]\to \P_2(\Om)$ is a weak solution of \eqref{eqn:agg-eqn2} if
	\[
		P_x(-\nabla V - \nabla K * \mu) \in L^1_{\loc}([0,T]; L^2(\mu(t))),
	\]
and for all $\varphi \in C_c^\infty(\Om \times (0,T))$,
	\[
		\int_0^T \! \int_{\Om} \big(\pt_t \varphi + \nabla \varphi \cdot P_x(-\nabla V - \nabla K * \mu) \big) \,d\mu \, dt=0, \qquad\quad \mu(0)=\mu^0,
	\]
where $\supp \mu(t) \subset \Om$ for all $t\in[0,T]$.
\end{defn}

\subsection{Numerical calculation of the $2$-Wasserstein distance}
\label{subsec:NumCalcW2}

In one dimension, the 2-Wasserstein distance between a classical density function and a sum of Dirac-deltas can be computed as follows. 

For $\rho_1,\,\rho_2 \in \P_2^a(\Om)$
the 2-Wasserstein distance between them can be written as
	\beqn \label{eqn:Wass_comp}
		d_W^2(\rho_1,\rho_2) = \inf_{\tbf} \int_{\Rd} |\tbf(x)-x|^2\rho_1(x)\,dx,
	\eeqn
where the infimum is taken over all maps $\tbf$ transporting $\rho_1$ to $\rho_2$. We say that a map $\tbf:\Om \to \Om$ transports $\rho_1$ into $\rho_2$ if
	\beqn \label{eqn:transp_map}
		\int_{x\in A} \rho_2(x)\,dx = \int_{\tbf(x)\in A} \rho_1(x)\,dx \qquad \text{for all bounded subsets } A \text{ of }\Rd.
	\eeqn
Then by a theorem due to Brenier (cf. \cite[Theorem 2.12]{Villani}) there exists a unique optimal transport map $\tbf_{\rho_1}^{\rho_2}$ which attains the infimum above, and it can be written as $\tbf_{\rho_1}^{\rho_2} = \nabla \psi$ for some convex function $\psi$. It is noted in \cite{BeBr2000} that in one-dimension the map $\tbf_{\rho_1}^{\rho_2}$ is non-decreasing by the convexity of $\psi$, and $\tbf_{\rho_1}^{\rho_2}$ can be determined entirely by the condition	
	\[
		\int_{x<\tbf_{\rho_1}^{\rho_2}(y)} \rho_2(x)\,dx = \int_{x<y} \rho_1(x)\,dx.
	\]

The considerations above also hold when one of the densities is a delta-measure. To compute the 2-Wasserstein distance between a density function and a sum of Dirac-deltas in one-dimension we will proceed as in \cite{BeBr2000}, and compute the optimal transport map $\tbf_{\rho}^{\mu}$ from a probability density $\rho$, compactly supported on $[0,L]$, into $\mu=\sum_{i=1}^n s_i \delta_{y_i}$. Since $\tbf_{\rho}^{\mu}$ is a non-decreasing function, using \eqref{eqn:transp_map}, we can find the optimal transport map by finding a partition $0=x_0<x_1<\cdots<x_n=L$ of the interval $[0,L]$, where the weight of the Dirac-mass at $y_i$ is obtained by integrating the density $\rho$ over the subinterval $[x_{i-1},x_i]$:
	\[
		s_i = \int_{[x_{i-1},x_i]} \rho(x)\,dx \qquad \text{for all } i=1,\ldots,n.
	\] 
Then we define $y_i=\tbf_{\rho}^{\mu}(x)$ for $x\in [x_{i-1},x_i]$ for all $i=1,\ldots,n$. Returning to \eqref{eqn:Wass_comp} we see that
	\[
		d_W^2(\rho,\mu) = \int_{[0,L]} |\tbf_{\rho}^{\mu}(x)-x|^2 \rho(x)\,dx = \sum_{i=1}^n \int_{[x_{i-1},x_i]} |y_i - x|^2 \rho(x)\,dx.
	\]


\section{Zero-diffusion limit: Analysis} \label{sec:analysis}

In this section we present the proofs of Theorems \ref{thm:conv} and \ref{thm:conv_min}. We start with the result regarding the convergence of weak solutions.

\subsection{Proof of Theorem \ref{thm:conv}}
We first note that the existence and uniqueness of weak solutions to model \eqref{eqn:eps-eqn} follows from the general theory of gradient flows in spaces of probability measures (see Section 11.2 and in particular Theorem 11.2.8 in \cite{AGS2005}). Specifically, under the assumptions made in Theorem \ref{thm:conv}, there exists a unique weak solution to \eqref{eqn:eps-eqn} which is a gradient flow of energy $\E_\nu$ \cite[Proposition 5.3]{CaCrPa17}. Also, existence and uniqueness of weak solutions to the aggregation equation \eqref{eqn:agg-eqn2} on domains with boundaries was proved in \cite[Theorems 1.5 and 1.6]{CarrilloSlepcevWu2016}), also using gradient flows framework and theory. 

Our proof of Theorem \ref{thm:conv} is adapted from the proof of \cite[Theorem 4.3]{Zh2017}, where linear diffusion is considered instead;
with the inclusion of a \emph{nonlinear} diffusion term in our case, some nontrivial modifications are needed as presented below.

\begin{proof}[Proof of Theorem \ref{thm:conv}]
For $0<\nu<1$, let $\mu_\nu(t)\in\Pa(\Om)$ be a weak solution of \eqref{eqn:eps-eqn} with density $\rho_\nu(t)$, i.e., $\mu_\nu(t) = \rho_\nu(t)\,dx$, and let $\mu(t)\in\P_2(\Om)$ be a weak solution of \eqref{eqn:agg-eqn2} for all $t\in[0,T]$. 

Consider a vector $\xibf_\nu\in L^2(\mu_\nu;\Om)$ in the Fr\'{e}chet subdifferential $\pt\E_\nu(\mu_\nu)$ of $\E_\nu$ at $\mu_\nu$. By the $\lambda$-convexity of $\E_\nu$ and the characterization of the subdifferential of $\E_\nu$ via an HWI-type inequality (see \cite[Section 10.1.1]{AGS2005} and equation (10.1.7), in particular), one has
	\beqn \label{eqn:subdiff-1}
		\E_\nu(\sigma) - \E_\nu(\mu_\nu) \geq \int_{\Om} \xibf_\nu \cdot (\tbf_{\mu_\nu}^\sigma - \ibf)\,d\mu_\nu  + \frac{\lambda}{2}d_W^2(\sigma,\mu_\nu),
	\eeqn
for any $\sigma\in \Pa(\Om)$.	
Similarly, consider a vector $\xibf$ in the Fr\'{e}chet subdifferential $\pt\E(\mu)$ of $\E$ at $\mu$ where $\E$ is defined by \eqref{eqn:energy}. Also by $\lambda$-convexity of $\E$, 
for any $\tilde{\sigma}\in\P_2(\Om)$,
	\beqn \label{eqn:subdiff-2}
		\E(\tilde{\sigma})-\E(\mu)  \geq \int_{\Om} \xibf \cdot (\tbf_{\mu}^{\tilde{\sigma}} - \ibf)\,d\mu + \frac{\lambda}{2}d_W^2(\tilde{\sigma},\mu).
	\eeqn

Now, using Lemma \ref{lem:modify} with base $\mu_\nu(t)$ and $s=\nu^{1/(d+2)}$ we obtain the existence of $\mu_s(t)\in\Pa(\Om)$ with $\mu_s(t) = \rho_s(t)\,dx$ so that for any $t\in [0,T]$,
	\[
		d_{\mu_\nu}(\mu_s(t),\mu(t)) \leq C_\Om \, \nu^{1/(d+2)} \qquad \text{and} \qquad \|\rho_s(t)\|_{L^\infty(\Om)} \leq \nu^{-d/(d+2)}.
	\]
We note here that the proof of \cite[Lemma 4.9]{Zh2017} shows explicitly that, although the base $\mu_\nu$ and the choice of $s$ depend on $\nu$, the constant $C_\Om$ in the estimate of $d_{\mu_\nu}$ is independent of $\nu>0$. 

Then, taking $\sigma=\mu_s$ in \eqref{eqn:subdiff-1} yields
		\begin{align}
			\E_\nu(\mu_s) - \E_\nu(\mu_\nu) &\geq \int_{\Om} \xibf_\nu \cdot (\tbf_{\mu_\nu}^{\mu_s} - \ibf)\,d\mu_\nu  + \frac{\lambda}{2}d_W^2(\mu_s,\mu_\nu) \nonumber \\
																	  &= \int_{\Om} \xibf_\nu \cdot  (\tbf_{\mu_\nu}^{\mu} - \ibf)\,d\mu_\nu
																	 + \int_{\Om} \xibf_\nu \cdot  (\tbf_{\mu_\nu}^{\mu_s}-\tbf_{\mu_\nu}^{\mu})\,d\mu_\nu  +  \frac{\lambda}{2}d_W^2(\mu_s,\mu_\nu) \nonumber \\
																	  & \geq \int_{\Om}\!\int_{\Om} \xibf_\nu(y) \cdot (x-y)\,d\gamma_\nu(x,y) + \int_{\Om} \xibf_\nu \cdot  (\tbf_{\mu_\nu}^{\mu_s}-\tbf_{\mu_\nu}^{\mu})\,d\mu_\nu \label{eqn:subdiff-3} \\
																	  &\quad    - C d_W^2(\mu,\mu_\nu) - \tilde{C} \nu^{2/(d+2)}, \nonumber
		\end{align}
where $\gamma_\nu$ denotes the optimal transport plan between $\mu$ and $\mu_\nu$, and $C, \tilde{C}$ denote two generic positive constants. The last inequality is immediate for $\lambda \geq 0$ and follows from triangle inequality for $\lambda<0$.

Returning to \eqref{eqn:subdiff-2}, we take $\tilde{\sigma}=\mu_\nu$ and get
		\beqn \label{eqn:subdiff-4}
			\E(\mu_\nu) - \E(\mu) \geq \int_{\Om} \! \int_{\Om} \xibf(x)\cdot (y-x)\,d\gamma_\nu(x,y) + \frac{\lambda}{2} d_W^2(\mu,\mu_\nu).
		\eeqn
	
By H\"{o}lder inequality with respect to the measure $\mu_\nu$,
		\beqn \label{eqn:Holder0}
			\left| \int_{\Om} \xibf_\nu \cdot (\tbf_{\mu_\nu}^{\mu_s}-\tbf_{\mu_\nu}^{\mu})\,d\mu_\nu \right| \leq \left( \int_{\Om}|\xibf_\nu|^2\,d\mu_\nu \right)^{1/2} \left(\int_{\Om} |\tbf_{\mu_\nu}^{\mu_s}-\tbf_{\mu_\nu}^{\mu}|^2\, d\mu_\nu \right)^{1/2}.
		\eeqn
By \cite[Theorem 11.3.2]{AGS2005}, $\mu_\nu$ is a gradient flow solution; hence, we can choose 
	\[
		\xibf_\nu = - \mathbf{v}_\nu = \frac{\nu^\alpha \, m}{m-1} \nabla \rho_\nu^{m-1} + \nabla K * \mu_\nu + \nabla V,
	\]
and by \cite[Theorem 11.3.4]{AGS2005} we get that $\int_{\Om} |\xibf_\nu|^2\,d\mu_\nu $ is in $L^\infty_{\text{loc}}(0,\infty)$. Moreover, we can adapt \cite[Lemma 4.7]{Zh2017} to our case, where the idea is to control the time-step discretization of $\xibf_\nu$ using the JKO scheme, and obtain that $\int_{\Om} |\xibf_\nu|^2\,d\mu_\nu $ is uniformly bounded (in $\nu$) in $L^\infty_{\text{loc}}(0,\infty)$. Therefore, using the fact that
	\[
		\left(\int_{\Om} |\tbf_{\mu_\nu}^{\mu_s}-\tbf_{\mu_\nu}^{\mu}|^2\, d\mu_\nu \right)^{1/2} = d_{\mu_\nu}(\mu,\mu_s)
	\]
defines the pseudo-Wasserstein distance induced by $\mu_\nu$, \eqref{eqn:Holder0} yields
	\beqn \label{eqn:Holder}
		\left| \int_{\Om} \xibf_\nu \cdot (\tbf_{\mu_\nu}^{\mu_s}-\tbf_{\mu_\nu}^{\mu})\,d\mu_\nu  \right| \leq C d_{\mu_\nu}(\mu,\mu_s)  \leq C \nu^{1/(d+2)}.
	\eeqn
Also, by the HWI-type inequality for $\E$ (see e.g. \cite[Proposition 2.5]{Craig17}), and the fact that $d_W (\cdot,\cdot)\leq d_{\mu_\nu}(\cdot,\cdot)$, we have
	\beqn \label{eqn:HWI}
		\E(\mu_s) - \E(\mu) \leq C d_{\mu_\nu}(\mu_s,\mu) \leq C \nu^{1/(d+2)}.
	\eeqn
Combining \eqref{eqn:subdiff-3} and \eqref{eqn:subdiff-4}, and using the estimates \eqref{eqn:Holder} and \eqref{eqn:HWI}, we then get 
	\begin{multline} \nonumber
		\int_{\Om}\!\int_{\Om} (-\xibf(x)+\xibf_\nu(y)) \cdot (x-y)\,d\gamma_\nu(x,y) \leq \Se_\nu(\mu_s(t)) - \Se_\nu(\mu_\nu(t))+ C  d_W^2(\mu,\mu_\nu) + \tilde{C} \nu^{1/(d+2)},
	\end{multline}
where $\Se_\nu(\mu(t)) = \frac{\nu^\alpha}{m-1}\int_{\Om}\rho^m(t)\,dx$ denotes the nonlinear diffusion part in the energy functional $\E_\nu$. Since $\rho_s(t) \leq \nu^{-d/(d+2)}$ pointwise, and since $\Om$ is bounded by (H4), we have that
	\[
		\Se_\nu(\mu_s(t)) \leq C \nu^{\alpha - dm/(d+2)}.
	\]
On the other hand, since $\Se_\nu$ is nonnegative, $\Se_\nu(\mu_\nu(t)) \geq 0$, and we get
	\beqn	\label{eqn:xi-1}
		\int_{\Om}\!\int_{\Om} (-\xibf(x)+\xibf_\nu(y)) \cdot (x-y)\,d\gamma_\nu(x,y) \leq C d_W^2(\mu,\mu_\nu) + \tilde{C} \nu^\beta,
	\eeqn
where
	\[
		\beta \,:=\, \min \left\{ \alpha - \frac{dm}{d+2},\, \frac{1}{d+2} \right\}.
	\]
Inequality \eqref{eqn:xi-1} holds for any vector $\xibf \in \pt\E(\mu)$. Take now $\xibf = \nabla K \ast \mu + \nabla V$. An expansion of the energy $\E$ around $\mu$ shows that this choice of $\xibf$ satisfies \eqref{eqn:subdiff-1} (also see equation (10.1.7) in \cite{AGS2005}) as $\E$ is $\lambda$-convex; hence, we get that $\nabla K \ast \mu + \nabla V \in \pt\E(\mu)$.
	
Also, by \cite[Lemma 4.3.4 and Theorem 8.4.7]{AGS2005}, we have
	\[
		\frac{d}{dt} d_W^2(\mu,\mu_\nu) = \int_{\Om}\!\int_{\Om} (\mathbf{v}(x)-\mathbf{v}_\nu(y))\cdot(x-y)\,d\gamma_\nu(x,y),
	\]
where $\mathbf{v}$ and $\mathbf{v}_\nu$ are tangent vector fields of $\mu$ and $\mu_\nu$ respectively. In fact, $\mathbf{v}=P_x(-\xibf)$, and since $\Om$ is convex by (H4), for all $x\in\pt\Om$ and for all $y\in\Om$,
	\[
		\big(P_x(-\xibf)(x)+\xibf(x) \big) \cdot (x-y) \leq 0.
	\]  
Also, $\supp (P_x(-\xibf)+\xibf) \subset \pt\Om $. Therefore, we get
	\beqn \label{eqn:xi-2}
		\frac{d}{dt} d_W^2(\mu,\mu_\nu) \leq  \int_{\Om}\!\int_{\Om} (-\xibf(x) + \xibf_\nu(y)) \cdot (x-y)\,d\gamma_\nu(x,y).
	\eeqn
Now, \eqref{eqn:xi-1} and \eqref{eqn:xi-2}, along with $d_W^2(\mu(0),\mu_\nu(0))=0$, imply by Gr\"{o}nwall's inequality that for all $t\in[0,T]$:
	\[
		d_W^2(\mu(t),\mu_\nu(t))  \leq \tilde{C} \nu^\beta \,t\,e^{C t}
	\]
for some constants $C, \tilde{C}>0$ that are independent of $\nu>0$.
\end{proof}

\subsection{Proof of Theorem \ref{thm:conv_min}}

We start by making the following remark:
\begin{remark}
\label{rem:bddbelow_lsc}
Any function $K$ satisfying the hypotheses (H1) and (H2B) is bounded from below and lower semicontinuous. Consequently, without loss of generality, we will assume that $K$ is nonnegative, since the minimizers of an interaction energy $\E$ with potential $K$ are the same as the minimizers of $\E$ with potential $K - \inf K$.
\end{remark}

Before we prove Theorem \ref{thm:conv_min}, the first question that needs to be addressed is the existence of minimizers of energies $\E$ and $\E_\nu$. In particular, this becomes a nontrivial issue when $\Om$ is an unbounded set, as any minimizing sequence could potentially escape to infinity due to loss of compactness. Following arguments in \cite{ChFeTo2015,SiSlTo2015}, we will establish the existence of minimizers in $\P_2(\Om)$ via Lions' concentration compactness lemma and compactness of the support of minimizers. We state Lions' lemma here for readers' convenience, followed by the result on the existence of minimizers for energies $\E$ and $\E_\nu$.

\blemma[{Concentration compactness lemma for measures (cf. \cite{Lions84}, \cite[Section 4.3]{Struwe})}]\label{lem:conc_comp}
		Let $\{\mu_n\}_{n\in\mathbb{N}}$ be a sequence of probability measures on $\Rd$. Then there exists a subsequence $\{\mu_{n_k}\}_{k\in\mathbb{N}}$ satisfying one of the three following possibilities:
		\vspace{0.2cm}
	\begin{itemize}
		\item[(i)] \emph{(tightness up to translation)} There exists a sequence $\{y_k\}_{k\in\mathbb{N}}\subset\Rd$ such that for all $\e>0$ there exists $R>0$ with the property that
			\[
				\int_{B_R(y_k)}\,d\mu_{n_k}(x) \geq 1-\e \qquad {\hbox{\rm  for all $k$.}}
			\]
		 \item[(ii)] \emph{(vanishing)} $\disp \lim_{k\arrow\infty} \sup_{y\in\Rd} \int_{B_R(y)}\,d\mu_{n_k}(x)=0$, for all $R>0$; \vspace{0.2cm}
		 \item[(iii)] \emph{(dichotomy)} There exists $\alpha\in(0,1)$ such that for all $\e>0$, there exist a number $R>0$ and a sequence $\{x_k\}_{k\in\mathbb{N}}\subset\Rd$ with the following property:\\
		 
		 	\vspace{-0.2cm}
		 Given any $R^\pr>R$ there are nonnegative measures $\mu_{1,k}$ and $\mu_{2,k}$ such that
		 	\vspace{0.2cm}
			\begin{itemize}
		 		\item[] $0\leq \mu_{1,k} + \mu_{2,k} \leq \mu_{n_k}$\,, \quad $\supp(\mu_{1,k})\subset B_R(x_k)$,\quad $\supp(\mu_{2,k})\subset \Rd\setminus B_{R^\pr}(x_k)$,
				\vspace{0.2cm}
		 		\item[] $\disp \limsup_{k\arrow\infty} \left(\left|\alpha-\int_{\Rd}d\mu_{1,k}(x)\right|+\left|(1-\alpha)-\int_{\Rd}d\mu_{2,k}(x)\right|\right)\leq \e$.
		 	\end{itemize}
	\end{itemize}
\elemma


The existence of minimizers for energy functionals of type $\E_\nu$ was proved by Lions in his seminal paper \cite{Lions84} for nonnegative and decaying interaction potentials $K$ in some Marcinkiewicz (or weak $L^p$) space. Here we prove the existence of minimizers for a more general class of interaction potentials, which includes attractive-repulsive potentials that grow to infinity, and also establish the compactness of the support of minimizers.

\bthm[Existence of minimizers] \label{thm:exist}
	Assume $K$ and $V$ satisfy the hypotheses \emph{(H1), (H2B)}, and \emph{(H3B)}. For any $m>1$ and $\nu>0$ the energies $\E_\nu$ defined by  \eqref{eqn:reg_energy} admit a compactly supported minimizer in $\P_2(\Om)$. The same statement holds true for the energy $\E$ in \eqref{eqn:energy}.
\ethm

\begin{proof}
We will prove the theorem for regularized energies $\E_\nu$. The existence of minimizers of $\E$ follows by the same arguments since the additional term in $\E_\nu$ is positive and lower semicontinuous.

Let $\{\mu_n\}_{n\in\mathbb{N}} \subset \P(\Rd)$ be a minimizing sequence of $\E_\nu$, that is, 
	\[\lim_{n\to\infty} \E_\nu(\mu_n)=\inf_{\mu\in\P(\Rd)} \E_\nu(\mu).\]
We can assume that for $n$ sufficiently large $\mu_n$ is absolutely continuous with respect to the Lebesgue measure with density $\rho_n\in L^1\cap L^m$ since otherwise $\E_\nu(\mu_n) \to \infty$ as $n\to\infty$. Also note that for notational convenience we omit the dependence on $\nu$ of the minimizing sequence $\mu_n$. Nevertheless, the dependence on $\nu$ will be explicitly indicated on $\mu_\nu$,  the $n \to \infty$  limit (on a subsequence) of $\mu_n$. 
\medskip

If $\Om$ is bounded, the conclusion is immediate. Indeed, $\E_\nu(\mu_n)$ is bounded above (for large $n$), and by the assumptions that $V$ is bounded from below and $K$ is positive, we get
	\beqn \label{eqn:est-Lm}
		\frac{\nu^\alpha}{m-1}\int_{\Om} \rho_n^m\,dx \leq - \frac{1}{2} \int_{\Om}\!\int_{\Om} K(x-y)\,\rho_n(x)\rho_n(y)\,dxdy - \int_{\Om} V(x) \rho_n(x)\,dx + C \leq \tilde{C},
	\eeqn
where $C, \tilde{C}$ are positive constants that do not depend on $n$. Hence, $\|\rho_n\|_{L^m(\Om)}$ are uniformly bounded for $n$ sufficiently large, and therefore, using boundedness of $\Om$ we get a subsequence of $\rho_n$ which converges weakly to a function $\rho_\nu$ in $L^m(\Om)$. This implies that $\rho_n \warrow \rho_\nu$ also in $\P(\Om)$ with respect to weak-* topology. By convexity of the first term, Remark \ref{rem:bddbelow_lsc}, and the Portmanteau Theorem (cf. \cite[Theorem 1.3.4]{van1996weak}) the first two terms of the energy $\E_\nu$ are lower semicontinuous with respect to weak-* convergence in $\P(\Om)$. The last term is lower semicontinuous again by the Portmanteau theorem since $V$ satisfies (H3B). Therefore $\mu_\nu = \rho_\nu\,dx \in \P(\Om)$ minimizes $\E_\nu$.

\medskip
Now, suppose $\Om$ is unbounded, and suppose $\{\rho_n\}_{n\in\mathbb{N}}$ has a subsequence which ``vanishes''. Since that subsequence is also a minimizing sequence we can assume that $\{\rho_n\}_{n\in\mathbb{N}}$ vanishes.
  Then  for any $\delta>0$ and for any $R>0$ there exists $N\in\mathbb{N}$ such that for all $n>N$ and for all $x\in\Om$
 	\[
 		\int_{\Om\setminus B_R(x)}\rho_n\,dx \geq 1-\delta.
 	\]
This implies that for $n>N$,
	\[
		\iint_{|x-y|\geq R} \rho_n(x)\rho_n(y)\,dxdy = \int_{\Om}\!\left(\int_{\Om\setminus B_R(x)} \rho_n(y)\,dy\right)\rho_n(x)\,dx \geq 1-\delta.
	\]
Given $M \in \R$, by condition (H2B) there exists $R>0$ such that for all $r \geq R$, $K(r) \geq M$. Let $\delta \in (0, \frac12)$, and take $N$ corresponding to $\delta$ and $R$.
Since $K \geq 0$ by Remark \ref{rem:bddbelow_lsc}, and $V \geq -C_V$ for some $C_V>0$ by (H3B), one has
	\beqns
			\begin{aligned}
				\E_\nu(\rho_n) &\geq  \frac{1}{2} \iint_{|x-y| < R} K(x-y)\rho_n(x)\rho_n(y)\,dxdy + \frac{1}{2} \iint_{|x-y|\geq R} K(x-y)\rho_n(x)\rho_n(y)\,dxdy - C_V \\
										 &\geq \frac{1}{2}(1-\delta)M - C_V,
			\end{aligned}
	\eeqns
for all $n>N$. Letting $M \to \infty$ implies $\E_\nu(\rho_n)\arrow\infty$. This contradicts the fact that $\rho_n$ is a subsequence of a minimizing sequence of $\E_\nu$. Thus, ``vanishing'' does not occur.

Next we show that ``dichotomy'' is also not an option for a minimizing sequence. Suppose,  that ``dichotomy'' occurs.  As before we can assume that the subsequence along which dichotomy occurs is the whole sequence. Let $\e>0$ be fixed, and let $R$, the sequence $\{x_k\}_{k\in\mathbb{N}}$ and measures
	\[
		\rho_{1,k}+\rho_{2,k} \leq \rho_{k}.
	\]
be as defined in Lemma \ref{lem:conc_comp}(ii).
For any $R'>R$, using Remark \ref{rem:bddbelow_lsc}, we obtain
	\beqns
		\begin{aligned}
		\liminf_{k\arrow\infty} \E(\rho_{n_k}) & \geq  \liminf_{k\arrow\infty} \frac{1}{2} \int_{B_R(x_{n_k})}\!\int_{B^c_{R^\pr}(x_{n_k})} K(x-y)\rho_{2,k}(x)\rho_{1,k}(y)\,dxdy - C_V \\
								 &\geq  \frac{1}{2} \inf_{r \geq R'-R}  K(r)\, (\alpha-\e)(1-\alpha-\e) - C_V ,
		\end{aligned}
	\eeqns
where $B^c_{R^\pr}(x_{n_k})$ simply denotes $\Om\setminus B_{R^\pr}(x_{n_k})$.

By (H3B), letting $R^\pr\arrow\infty$ yields that
	\[
		\liminf_{k\arrow\infty} \E(\rho_{n_k}) \geq \infty,
	\]
which contradicts the fact that $\{\rho_n\}$ is an energy minimizing sequence.

Therefore ``tightness up to translation'' is the only possibility. Hence there exists $y_k\in\Om$ such that for all $\e>0$ there exists $R>0$ with the property that
			\[
				\int_{B(y_k,R)} \rho_{n_k}\,dx \geq 1-\e \qquad \text{ for all } k.
			\]
Let
	\[
		\tilde{\rho}_{n_k}:=\rho_{n_k}(\cdot + y_k).
	\]
Then the sequence $\{\tilde{\rho}_{n_k}\}_{k\in\mathbb{N}}$ is tight. Therefore by Prokhorov's theorem (cf. \cite[Theorem 4.1]{Billingsley}) there exists a further subsequence of $\{\tilde{\rho}_{n_k}\}_{k\in\mathbb{N}}$ which we still index by $k$, and a measure $\mu_\nu \in\Prob(\Om)$ such that
	\[
		\tilde{\rho}_{n_k} {\rightharpoonup}\mu_\nu
	\]
in the weak-* topology of $\mathcal{P}(\Rd)$ as $k\arrow\infty$. On the other hand, note that since $V$ is either translation invariant or strictly increasing in each component, we see that either $\E_\nu(\tilde{\rho}_{n_k})=\E_\nu(\rho_{n_k})$ or $y_k\to 0$ as $k\to\infty$ for a minimizing sequence, and $\E_\nu(\tilde{\rho}_{n_k})\geq \E_\nu(\rho_{n_k})$ for sufficiently large $k$; hence, we can omit translations and conclude that $\rho_{n_k} \rightharpoonup \mu_\nu$ as $k\to\infty$. As in the bounded $\Om$ case, the energy $\E_\nu$ is lower semicontinuous by Portmanteau theorem; hence, $\mu_\nu \in\P(\Om)$ minimizes $\E_\nu$.


\medskip

Next, we will show that $\mu_\nu$ is in fact compactly supported; in particular, $\mu_\nu \in \P_2(\Om)$. From the definition of $\E_\nu$ it is clear that if $\mu_\nu$ minimizes $\E_\nu$ over $\P(\Om)$ then $\mu_\nu$ is absolutely continuous with respect to the Lebesgue measure with density $\rho_\nu \in L^1 \cap L^m$.

A simple first variation calculation shows that $\rho_\nu$ satisfies
	\beqn \label{eqn:first-var}
		\frac{m\nu^\alpha}{m-1} \rho_\nu^{m-1}(x) + \int_{\Om} K(x-y)\rho_\nu(y)\,dy + V(x) = \lambda,
	\eeqn
for some $\lambda\in\R$ and for all $x\in\supp\rho_\nu$. By the positivity of $\rho_\nu$, and since $V \geq -C_V$ by (H3B), we get
	\[
		\int_{\Om} K(x-y)\rho_\nu(y)\,dy  \leq \lambda + C_V,
	\]
for all $x\in\supp\rho_\nu$. On the other hand, for $|x|$ large,
	\[
		\int_{\Om} K(x-y)\rho_\nu(y)\,dy \geq \int_{\Om \cap \{ |y|< R\}} K(x-y)\rho_\nu(y)\,dy \geq C_R \inf\big\{K(z) \colon |z| \geq |x|-R\big\},
	\]
where $R>0$ is chosen large enough so that $C_R := \int_{\Om \cap \{|y| < R\}}\rho_\nu(y)\,dy>0$. Thus we have that $\lim_{|x|\to\infty} \int_{\Om}K(x-y)\rho_\nu(y)\,dy \to \infty$; hence, $\big\{x\in\Om \colon \int_{\Om} K(x-y)\rho_\nu(y)\,dy \leq \lambda+C_V \big\}$ is bounded. We conclude that $\mu_\nu$ has compact support, and also, bounded second moment.
\end{proof}

Next we prove the $\Gamma$-convergence of the energies $\E_\nu$ to the energy $\E$ as stated in the next theorem.

\bthm[$\Gamma$-convergence of energies]\label{thm:gamma_conv} 
Assume $K$ and $V$ satisfy the hypotheses \emph{(H1), (H2B)}, and \emph{(H3B)}. Also assume $\alpha > dm$. Let $\E_\nu$ and $\E$ be defined by \eqref{eqn:reg_energy} and \eqref{eqn:energy}, respectively.
\begin{itemize}
		\item[(i)] \emph{(Lower bound.)} For any $\{\mu_\nu\}_{\nu>0}\subset\P_2(\Om)$ and $\mu \in \P_2(\Om)$ with $\lim_{\nu\to 0}d_W(\mu_\nu,\mu)=0$ we have that
	\[
		\liminf_{\nu\to 0} \E_\nu(\mu_\nu) \geq \E(\mu).
	\]
		\item[(ii)] \emph{(Upper bound.)} For any $\mu\in\P_2(\Om)$ there exists  $\{\sigma_\nu\}_{\nu>0}\subset\P_2(\Om)$ such that
	\[
		\lim_{\nu\to 0}d_W(\sigma_\nu,\mu)=0 \quad\text{ and }\quad \lim_{\nu\to 0} \E_\nu(\sigma_\nu) = \E(\mu).
	\]
	\end{itemize}
\ethm

\begin{proof}
Let $\{\mu_\nu\}_{\nu>0}\subset \P_2(\Om)$ be a sequence such that $\lim_{\nu\to 0}d_W(\mu_\nu,\mu)=0$ for some $\mu\in\P_2(\Om)$. We can assume that $d\mu_\nu = \rho_\nu\,dx$ for $\rho_\nu\in L^1\cap L^m$ since otherwise the conclusion of part (i) trivially holds. Then
	\beqns
	\begin{aligned}
		\liminf_{\nu\to 0}\E_\nu(\rho_\nu) &\geq \liminf_{\nu\to 0} \frac{1}{2}\int_{\Om}\!\int_{\Om} K(x-y)\rho_\nu(x)\rho_\nu(y)\,dxdy + \int_{\Om} V(x)\rho_\nu(x)\,dx \\
																		& \geq \frac{1}{2} \int_{\Om}\!\int_{\Om} K(x-y)\,d\mu(x) d\mu(y) + \int_{\Om} V(x)\,d\mu(x) \\
																		& = \E(\mu)
	\end{aligned}
	\eeqns
by weak-* convergence of $\rho_\nu$ to $\mu$ and the Portmanteau theorem since $K$ is lower semicontinuous and bounded from below by Remark \ref{rem:bddbelow_lsc}, and $V$ is continuous and bounded from below by (H3B). Hence, part (i) follows.

Now, let $\mu\in\P_2(\Rd)$ be an arbitrary measure. Let $\varphi$ be any smooth function on $\Om$ such that $\int_{\Om} |x|^2\varphi(x)\,dx < \infty$, and define $\varphi_\nu = \nu^{-d}\varphi(x/\nu)$. Further define
	\[
		\sigma_\nu \,:=\, \varphi_\nu * \mu = \int_{\Om} \varphi_\nu (x-y)\,d\mu(y)
	\]
where the convolution operator $*$ is defined via an integral over $\Om$. Then, by \cite[Lemma 7.1.10]{AGS2005},
	\[
		d_W(\varphi_\nu*\mu,\mu) \to 0
	\]
as $\nu\to 0$. This also implies that $\varphi_\nu*\mu \rhu \mu$ in the weak-* topology of $\P(\Om)$. Also, note that $\sigma_\nu$ is in fact absolutely continuous with respect to the Lebesgue measure with density $\psi_\nu \in L^1 \cap L^m$ (cf. \cite[Remark 2.2]{CraigTopaloglu}).

Now, recall by (H2B) the pairwise interaction potential $K$ can be written as the sum of $K^r + K^a$ where $K^a$ is continuous and has at most quadratic growth. Then, using \cite[Proposition 2.6]{CraigTopaloglu} and the at most quadratic growth of $K^a$ and $V$, we get
	\beqns
		\begin{aligned}
				\limsup_{\nu\to 0} \E_\nu(\sigma_\nu) &\leq \limsup_{\nu\to 0} \Bigg( \frac{\nu^\alpha}{m-1}\int_{\Om} \psi_\nu^m(x)\,dx \\
																				 & \qquad\qquad +\frac{1}{2} \int_{\Om}\!\int_{\Om} (K^r+K^a)(x-y)\psi_\nu(x)\psi_\nu(y)\,dxdy + \int_{\Om} V(x)\psi_\nu(x)\,dx \Bigg) \\
																				 &\leq \limsup_{\nu\to 0} \Bigg( \frac{\nu^\alpha}{m-1}\int_{\Om} \psi_\nu^m(x)\,dx  \Bigg) \\
																				 & \qquad\qquad +\frac{1}{2} \int_{\Om} \! \int_{\Om} K(x-y)\,d\mu(x)d\mu(y) + \int_{\Om} V(x)\,d\mu(x).
		\end{aligned}
	\eeqns
Since $\|\varphi_\nu \|_{L^\infty(\Om)} \leq \nu^{-d}$ by definition and $\alpha > dm$ by assumption, by Young's inequality we get that the first term on the right-hand side of the inequality above satisfies
	\beqns
		\begin{aligned}
		\limsup_{\nu\to 0} \frac{\nu^\alpha}{m-1}\int_{\Om} \psi_\nu^m(x)\,dx &= \limsup_{\nu\to 0} \frac{\nu^\alpha}{m-1}\int_{\Om} \big(\varphi_\nu*\mu\big)^m(x)\,dx \\
																																	& \leq \limsup_{\nu\to 0} \frac{\nu^\alpha}{m-1} \nu^{-dm} \mu(\Om) = 0.
		\end{aligned}
	\eeqns
Hence, together with part (i),  we infer that $\lim_{\nu\to 0} \E_\nu(\sigma_\nu) = \E(\mu)$.
\end{proof}

The $\Gamma$-convergence of functionals provides the necessary structure so that one obtains the convergence of respective minimizers provided the functionals have sufficient compactness properties for their sequences of minimizers. We state this result in the next lemma. Proof of this lemma follows via Lions' concentration compactness lemma and by adapting the arguments in Theorem \ref{thm:exist}.

\blemma \label{lem:cmpct_min_seq}
Let $\{\mu_\nu\}_{\nu>0} \subset \P(\Rd)$ be a sequence such that for all sufficiently small $\nu>0$ the energies $\E_\nu(\mu_\nu) \leq C$ for some constant $C>0$. Then a subsequence of $\{\mu_\nu\}_{\nu>0}$ converges to a measure $\mu\in \P(\Rd)$ with respect to the weak-* topology.
\elemma

\medskip

As a result, we obtain the convergence of minimizers.

\begin{proof}[Proof of Theorem \ref{thm:conv_min}]
Let $\{\mu_\nu\}_{\nu>0}\subset\P_2(\Rd)$ be a sequence of minimizers of $\E_\nu$. Then there exists $C>0$ such that $\E_\nu(\mu_\nu) \leq C$ for $\nu>0$ sufficiently small. Hence, by Lemma \ref{lem:cmpct_min_seq}, there exists $\mu\in\P(\Rd)$ such that, up to a subsequence, $\mu_\nu \to \mu$ as $\nu\to 0$ in the weak-* topology of $\P(\Rd)$.

Now let $\sigma\in \P_2(\Rd)$ be arbitrary. By Theorem \ref{thm:gamma_conv}(ii), there exists a sequence $\{\sigma_\nu\}_{\nu>0}\subset \P_2(\Rd)$ and we have that
	\[
		\E(\mu) \leq \liminf_{\nu\to 0} \E_\nu(\mu_\nu) \leq \liminf_{\nu\to 0} \E_\nu(\sigma_\nu) = \E(\sigma).
	\]
However, since minimizers of $\E$ are compactly supported we have that $\inf_{\P(\Om)} \E = \inf_{\P_2(\Om)} \E$. Therefore $\mu$ minimizes $\E$ over $\P(\Om)$, and since it is compactly supported $\mu\in \P_2(\Om)$.

Lastly, we will show that, in fact, $\mu_\nu \to \mu$ in $\P_2(\Om)$. Note that 
	\[
		\supp \mu_\nu \subset \left\{ x\in\Om \colon \int_{\Om} K(x-y)\rho_\nu(y)\,dy \leq \lambda_\nu + C_V \right\}
	\]
for some $\lambda_\nu \in \R$ by \eqref{eqn:first-var}. On the other hand, for $\nu$ small enough we have, by assumption, that $\E_\nu(\mu_\nu) \leq C$. Hence, $\nu^\alpha\int_{\Om} \rho_\nu^m(x)\,dx \leq \tilde{C}(m-1)$, where $\tilde{C} = C + C_V$ (see \eqref{eqn:est-Lm}). Similarly, $\frac{1}{2} \int_{\Om}\!\int_{\Om} K(x-y)\rho_\nu(x)\rho_\nu(y)\,dxdy \leq \tilde{C}$.

Now multiplying both sides of \eqref{eqn:first-var} by $\rho_\nu$ and integrating over $\Om$ yields an upper bound on $\lambda_\nu$ independent of $\nu$:
	\[
		\lambda_\nu = \E_\nu(\mu_\nu) + \frac{1}{2} \int_{\Om}\!\int_{\Om} K(x-y)\rho_\nu(x)\rho_\nu(y)\,dxdy  + \nu^\alpha \int_{\Om} \rho_\nu^m(x)\,dx \leq C + \tilde{C}+ \tilde{C}(m-1).
	\]
Therefore the set $\left\{ x\in\Om \colon \int_{\Om} K(x-y)\rho_\nu(y)\,dy \leq \lambda_\nu + C_V \right\}$ is bounded uniformly in $\nu$; hence, for all sufficiently small $\nu>0$, $\supp \mu_\nu \subset \{x\in\Om \colon |x| \leq R\}$ for some large $R>0$.

\end{proof}

\section{Zero-diffusion limit: Numerics}
\label{sect:numerics}

The goals of the numerical studies presented in this section are the following. First, we bring numerical support to Theorems \ref{thm:conv} and \ref{thm:conv_min} in one dimension. Second, we demonstrate the desirable property of the diffusive model \eqref{eqn:eps-eqn} to bypass the unstable equilibria of the plain aggregation model \eqref{eqn:agg-eqn2}. The latter is one of the main motivations for this paper, that is, to rectify a flaw of the plain aggregation model in domains with boundaries. We also show similar numerical results with relaxed regularity assumptions on the interaction potential $K$, suggesting that hypothesis (H2A) might be too strong and the methods used in the proof of Theorem \ref{thm:conv} could potentially be refined to cover less regular potentials.

All numerical tests are in one dimension. The computational domain is the closed interval $[0,1.5]$. We do not consider any external potential, $V(x) \equiv 0$, and 
for the diffusive model \eqref{eqn:eps-eqn} we take $\alpha = 1$ and $m = 2$ throughout all of this section.  We will use the initial condition
\[\mu^0 = 4\hspace{0.04in}\One_{[0,0.25]}\]
for all simulations of models \eqref{eqn:eps-eqn} and \eqref{eqn:agg-eqn2}, where $\One_{[a,b]}$ denotes the characteristic function on $[a,b]$.

The plain aggregation model \eqref{eqn:agg-eqn2} is investigated using the particle method described in \cite{FeKo2017}. In \cite{CarrilloSlepcevWu2016}, the authors establish several important results concerning particle methods in domains with boundaries. In particular, they show the well-posedness of the approximating particle system and the weak convergence of its solutions to solutions of the PDE model \eqref{eqn:agg-eqn2}, as the number of particles approaches infinity. The method of enforcing the no-flux boundary condition is detailed in \cite{FeKo2017}. In one dimension as we have here, particles can move up to and onto the boundary point at the origin but cannot move past it. Once on the boundary, they can move off of it provided they have a positive velocity; otherwise they remain at origin.

The diffusive model \eqref{eqn:eps-eqn} is investigated using the finite volume method described in \cite{CaChHu2015}. The method preserves the non-negativity of solutions as well as the energy gradient flow structure, and it has been demonstrated to capture accurately the long-time behaviour and equilibria of model  \eqref{eqn:eps-eqn}. In particular, it works well for small diffusion values as in the present study, as it can deal robustly with metastable behaviour and large concentrations.   
To compare $2$-Wasserstein distance between solutions of the models \eqref{eqn:eps-eqn} and \eqref{eqn:agg-eqn2} we use the method detailed in Section \ref{subsec:NumCalcW2}. 

We consider two interaction potentials: a continuous (but non-differentiable) potential and a $C^2$ regularization of it. The continuous potential consists of Newtonian repulsion and quadratic attraction:
\begin{equation}
\label{eqn:C0intpot}
K(x) = \phi(x) + \frac{1}{2} x^{2},
\end{equation}
where $\phi(x)$ is the free-space Green's function for the negative Laplace
operator $-\Delta$ in one dimension:
\begin{equation}
\phi(x)=-\frac{1}{2}|x|.
\label{eqn:phi}
\end{equation}
Equilibria of model \eqref{eqn:agg-eqn2} with potential \eqref{eqn:C0intpot} in domains with boundaries was investigated in \cite{FeKo2017}. A key finding there was that, as $t \to \infty$, solutions to the plain aggregation model approach equilibria that are not minimizers of the energy.
Such unstable equilibria are disconnected, consisting of delta concentrations on the boundary and a free swarm in the interior of the domain. 


In order to satisfy the assumptions in Theorem \ref{thm:conv} (in particular assumption (H2A) on $K$), we also considered various smoothed versions of \eqref{eqn:C0intpot}.
We present below numerical results for a $C^2$ potential; similar results were also obtained using a $C^1$ regularization. Specifically, we apply regularization to $|x|$ by an even powered polynomial in $[-\epsilon,\epsilon]$ and get the $C^2$ potential:
\begin{equation}
\label{eqn:C2intpot}
K_\epsilon(x) = \phi_\epsilon(x) + \frac{1}{2} x^{2},
\end{equation}
where $\phi_\epsilon(x)$ is:
\begin{equation}\label{eqn:phi2}
\phi_\epsilon(x)=
\begin{cases}
\vspace{0.2cm}
\frac{1}{16\epsilon^3}x^4 - \frac{3}{8\epsilon}x^2 - \frac{3}{16}\epsilon & \quad x \leq \epsilon\\
-\frac{1}{2}|x| & \quad x > \epsilon.
\end{cases}
\end{equation}
In all simulations we set $\epsilon = 0.1$.  

In \cite{FeRa10} the authors investigated steady states of the plain aggregation model \eqref{eqn:agg-eqn2} with a $C^1$ regularized version of \eqref{eqn:C0intpot} and no external potential. They found that steady states were composed of a sum of delta masses, which could be stable or unstable depending on the mass distribution among the concentrations. This interaction potential is revisited in \cite{HuFe2013}, with the addition of diffusive terms. The scope of the investigations in \cite{HuFe2013} is much broader in fact, as the authors consider linear and nonlinear diffusion with power-law interaction potentials up to fourth order; in particular they study the diffusive model \eqref{eqn:eps-eqn} with potential \eqref{eqn:C0intpot}. Referencing the results in \cite{FeRa10}, it has been noted in \cite{HuFe2013} that when sufficiently small diffusion is added, steady states of multiple concentrations become multiple smoothed aggregates. Furthermore these states could possibly be continuous, piecewise smooth, and with compact support. 

In our simulations with the $C^2$ potential \eqref{eqn:C2intpot} we find again that solutions of the plain aggregation model approach asymptotically unstable steady states. In this case, the unstable equilibria consist of a sum of delta masses, with one such concentration at the boundary of the domain. This observation expands the conclusions from  \cite{FeKo2017} and suggests that the degeneracy of the plain aggregation model (with regard to evolution into unstable equilibria), is not related to the smoothness (or lack thereof) of the interaction potential, but it is rather generic. The plain aggregation model with smooth interaction potentials also needs rectification, which in this paper is provided by small nonlinear diffusion.



%
%
%

\subsection{Time evolution: $C^2$ interaction potential} \label{subsec:c2intpot}

We first present numerical support for Theorem \ref{thm:conv}, i.e., compare the $2$-Wasserstein distance between the solutions to \eqref{eqn:eps-eqn} and \eqref{eqn:agg-eqn2} and show that for fixed times, the distance decreases as $\nu$ decreases. Note that the estimate proved in Theorem \ref{thm:conv} is based on Gr\"{o}nwall's lemma, meaning that for a fixed $\nu$, the distance between the two solutions can potentially grow exponentially fast in time. For this reason the numerical check of Theorem \ref{thm:conv} is restricted to relatively early times, as large times would require simulations with diffusion values that are too small for numerical purposes. 

\begin{figure}[htb]
 \begin{center}
 \begin{tabular}{cc}
 \includegraphics[width=0.46\textwidth]{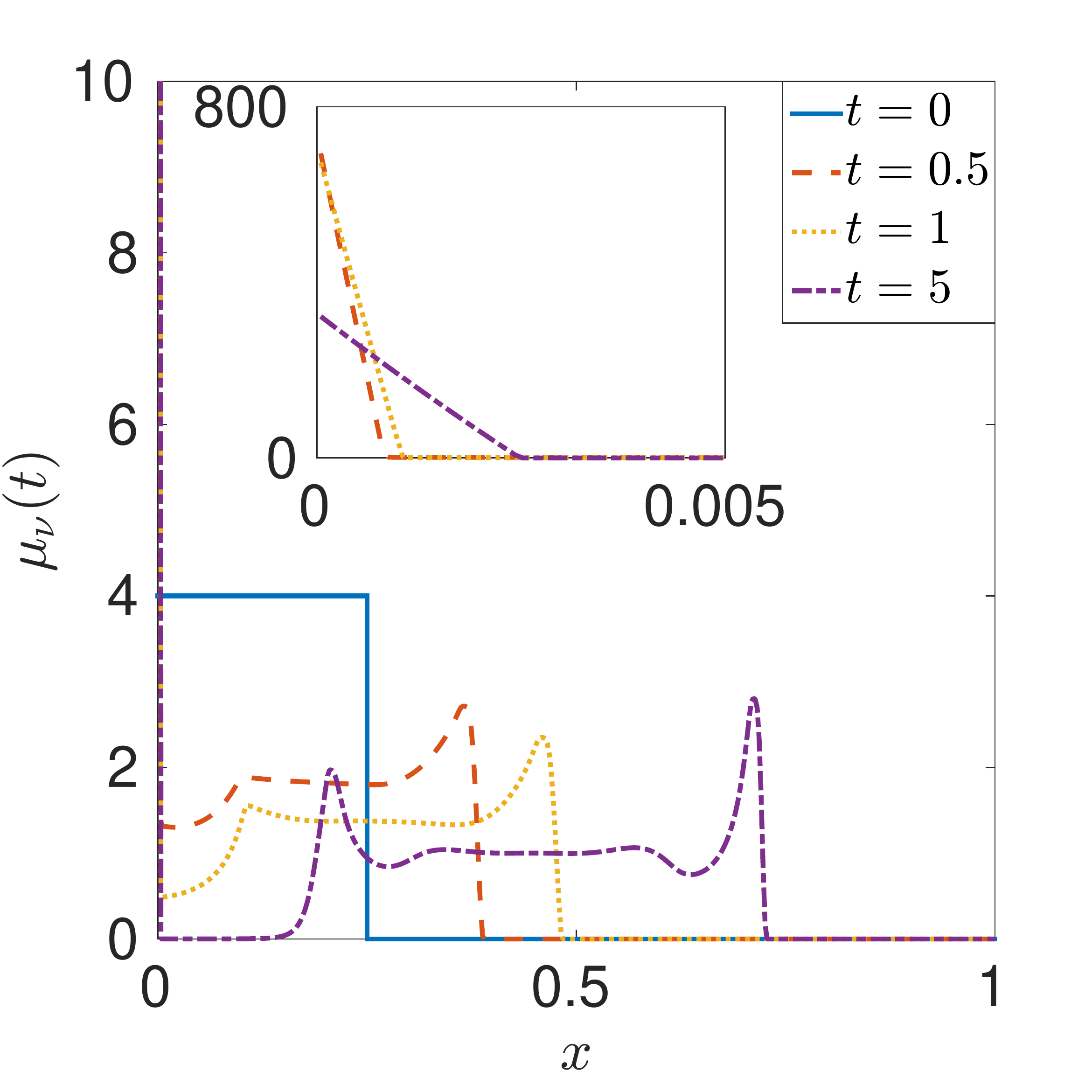} &
 \includegraphics[width=0.46\textwidth]{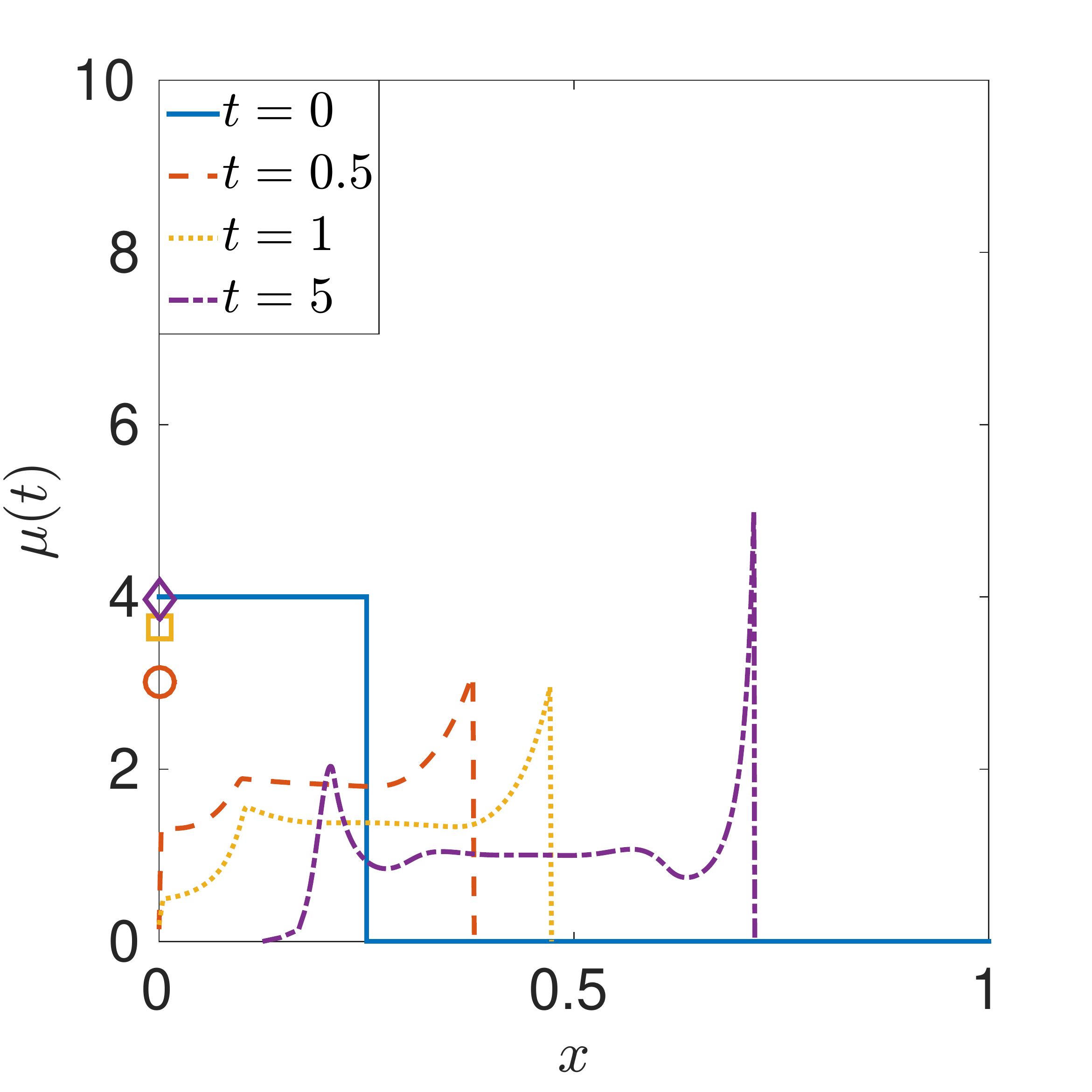} \\
 (a) & (b)
 \end{tabular}
\caption{Simulations with potential \eqref{eqn:C2intpot} showing early time dynamics. (a) Snapshots of the diffusive model \eqref{eqn:eps-eqn} with $\nu = 10^{-7}$. An insert has been included to show the layer of mass near the origin. (b) Snapshots of the plain aggregation model \eqref{eqn:agg-eqn2}. Concentrations at the origin are represented as circle, square, and diamond markers for $t = 0.5$, $t=1$, and $t=5$ respectively. The masses of concentrations have been magnified 10 times for clarity.}
\label{fig:SmoothstatesEarly}
\end{center}
\end{figure}

Indeed, for early times, the diffusive and plain aggregation models are qualitatively similar (see Figure~\ref{fig:SmoothstatesEarly}) and quantitatively, they remain close in the $2$-Wasserstein metric (see Table~\ref{tab:W2distEarly}). In both models we find that the initial mass begins separating, with some mass accumulating on, or near, the boundary and the rest moving away while remaining a single component. The notable difference is that the plain aggregation model forms delta concentrations of mass at the origin (Figure~\ref{fig:SmoothstatesEarly}(b)), whereas the diffusive model forms instead a thin, sharp layer of mass next to it (Figure~\ref{fig:SmoothstatesEarly}(a)). The latter is anticipated, as the measure-valued solutions of the diffusive model are absolutely continuous with respect to the Lebesgue measure.  

\begin{figure}[htb!]
  \begin{center}
 \begin{tabular}{cc}
 \includegraphics[width=0.46\textwidth]{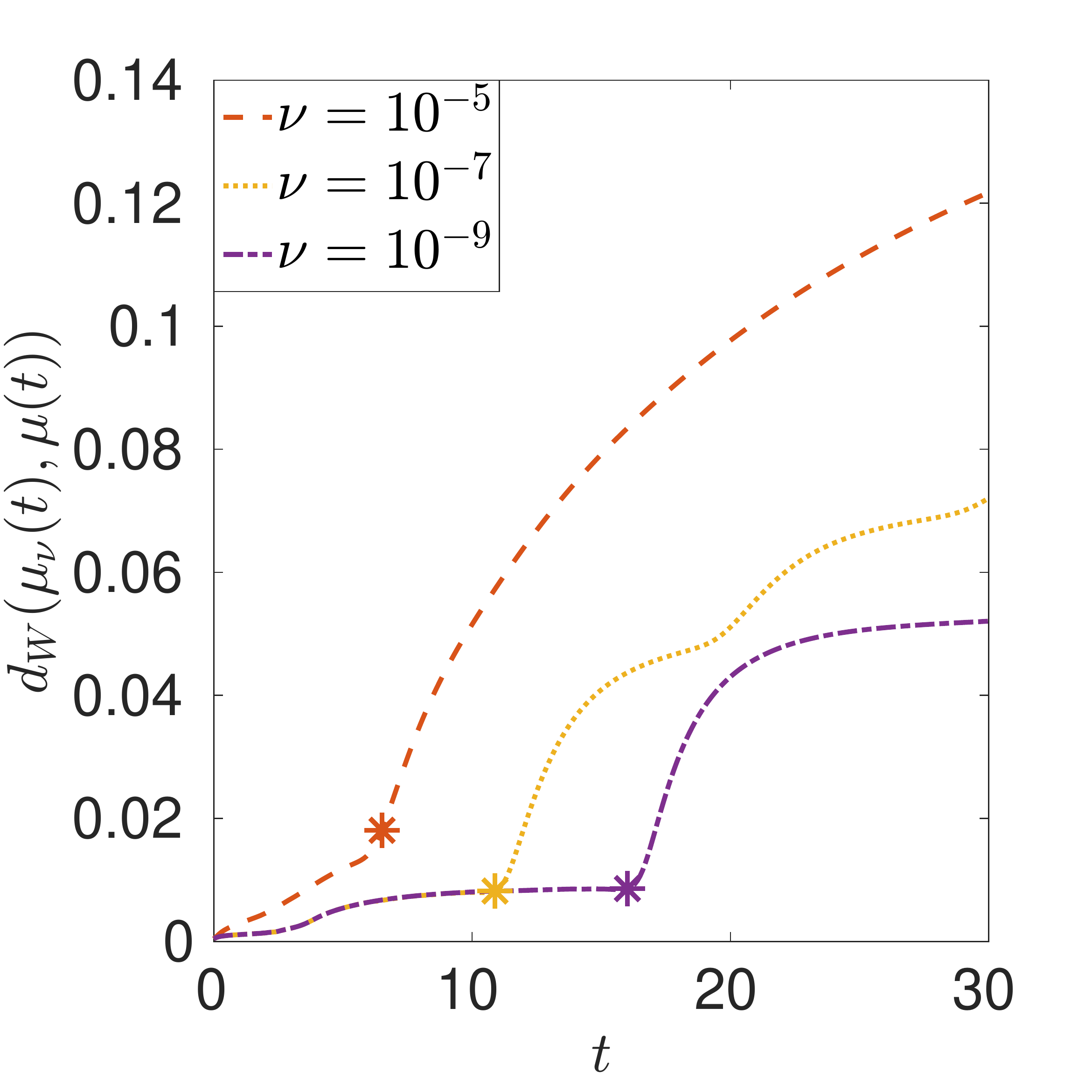} &
  \includegraphics[width=0.46\textwidth]{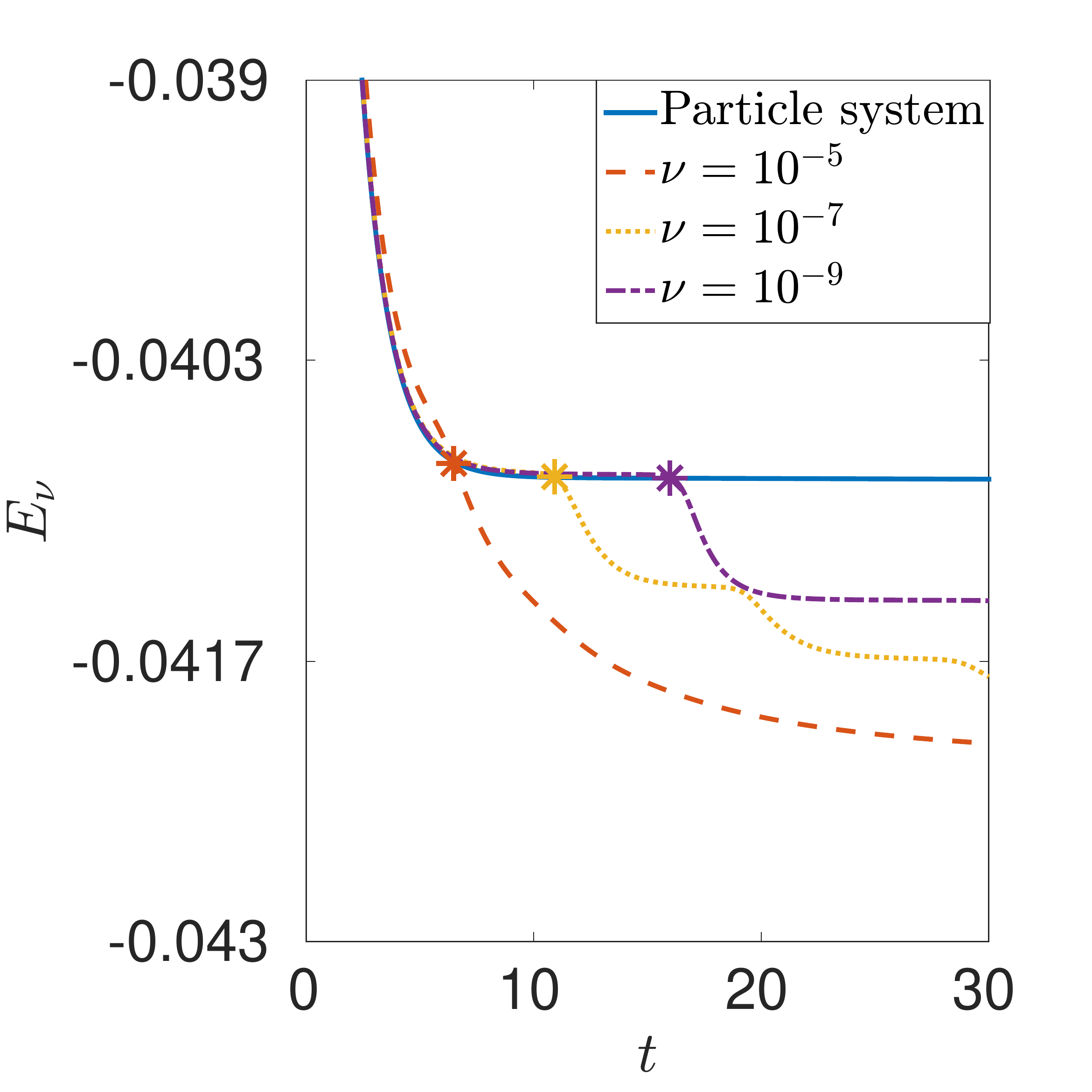} \\
 (a) & (b)
 \end{tabular} 
\caption{Results with potential \eqref{eqn:C2intpot}. (a) $2$-Wasserstein distance between the diffusive and plain aggregation solutions for various choices of $\nu$.  (b) Energy \eqref{eqn:reg_energy} of solutions to the diffusive model through time for various choices of $\nu$. Also included is the energy \eqref{eqn:energy} of the solution to the particle model through time (solid line). Star markers have been placed at $t = 6.5$, $t=10.9$ and $t=16$ for $\nu = 10^{-5}$, $\nu = 10^{-7}$ and $\nu = 10^{-9}$ respectively, corresponding to the times of the first mass transfer.} 
\label{fig:SmoothW2FVMvsP}
\end{center}
\end{figure}

Regarding the quantitative findings in Table~\ref{tab:W2distEarly}, we see two general trends.  First, the $2$-Wasserstein distance at a fixed time decreases as $\nu$ decreases towards zero, in support of Theorem~\ref{thm:conv}. Second, the $2$-Wasserstein distance between the two solutions grows as solutions are evolved through time; at early times they do so at a slow rate, for growth at later times see Figure~\ref{fig:SmoothW2FVMvsP}(a).

\bigskip

\begin{centering}
\begin{tabular}{|c|c|c|c|}
\hline 
$\nu$ & $t=0.5$ & $t=1$ & $t=5$ \\ 
\hline 
$10^{-3}$ & $2.1400\scie{-2}$ & $3.1896\scie{-2}$ & $8.6286\scie{-2}$
\\ 
\hline 
$10^{-4}$ & $7.1776\scie{-3}$ & $1.1110\scie{-2}$ & $4.7896\scie{-2}$
\\ 
\hline 
$10^{-5}$ & $3.8652\scie{-3}$ & $7.5722\scie{-3}$ & $3.3105\scie{-2}$
\\ 
\hline 
$10^{-6}$ & $3.4619\scie{-3}$ & $7.5057\scie{-3}$ & $3.3048\scie{-2}$
\\ 
\hline 
$10^{-7}$ & $3.4555\scie{-3}$ & $7.5043\scie{-3}$ & $3.2934\scie{-2}$
\\ 
\hline 
\end{tabular}\captionof{table}{$2$-Wasserstein distance $d_W(\mu_\nu(t),\mu(t))$ between solutions of the diffusive model and solutions of the plain aggregation model for various choices of $\nu$ and several early times.} \label{tab:W2distEarly} \par
\end{centering}
\medskip
 
The second major goal of these numerical simulations is to show that solutions of the diffusive model do not get trapped in equilibria that are not minimizers of the energy, as the plain aggregation model does. For this study we keep $\nu>0$ fixed and observe the long time evolution of solutions $\mu_\nu(t)$ of the diffusive model. As expected, we find that the distance between solutions $\mu_\nu(t)$ and $\mu(t)$ grows in time, with the caveat that the two solutions begin to differ substantially, both qualitatively and quantitatively, when mass near the origin in the diffusive model begins to move away from the boundary, into the interior of the domain. This mass transfer, a fundamental distinction between the two models,  will be highlighted throughout the discussion below.

\begin{figure}[thb!]
   \begin{center}
 \begin{tabular}{cc}
 \includegraphics[width=0.46\textwidth]{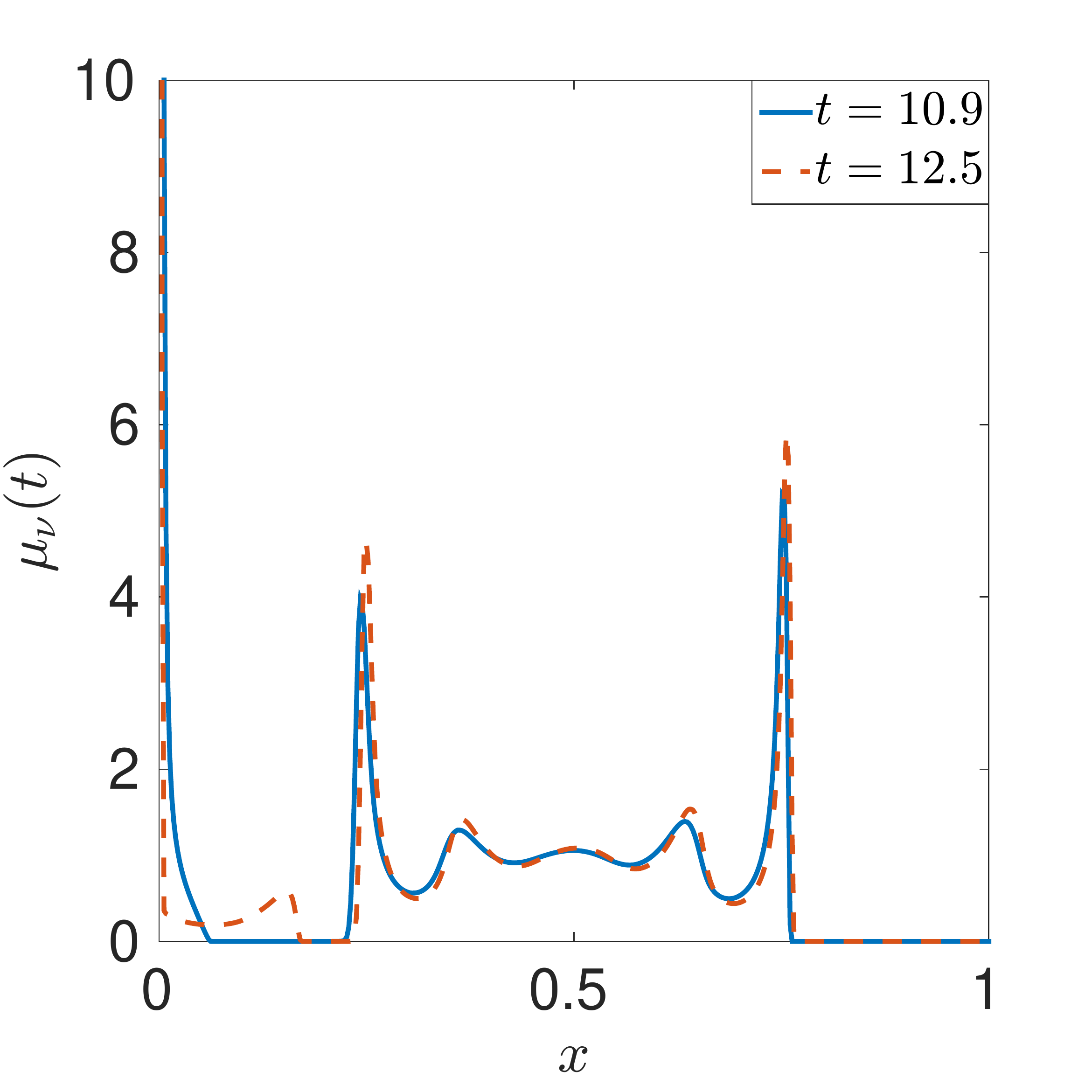} &
 \includegraphics[width=0.46\textwidth]{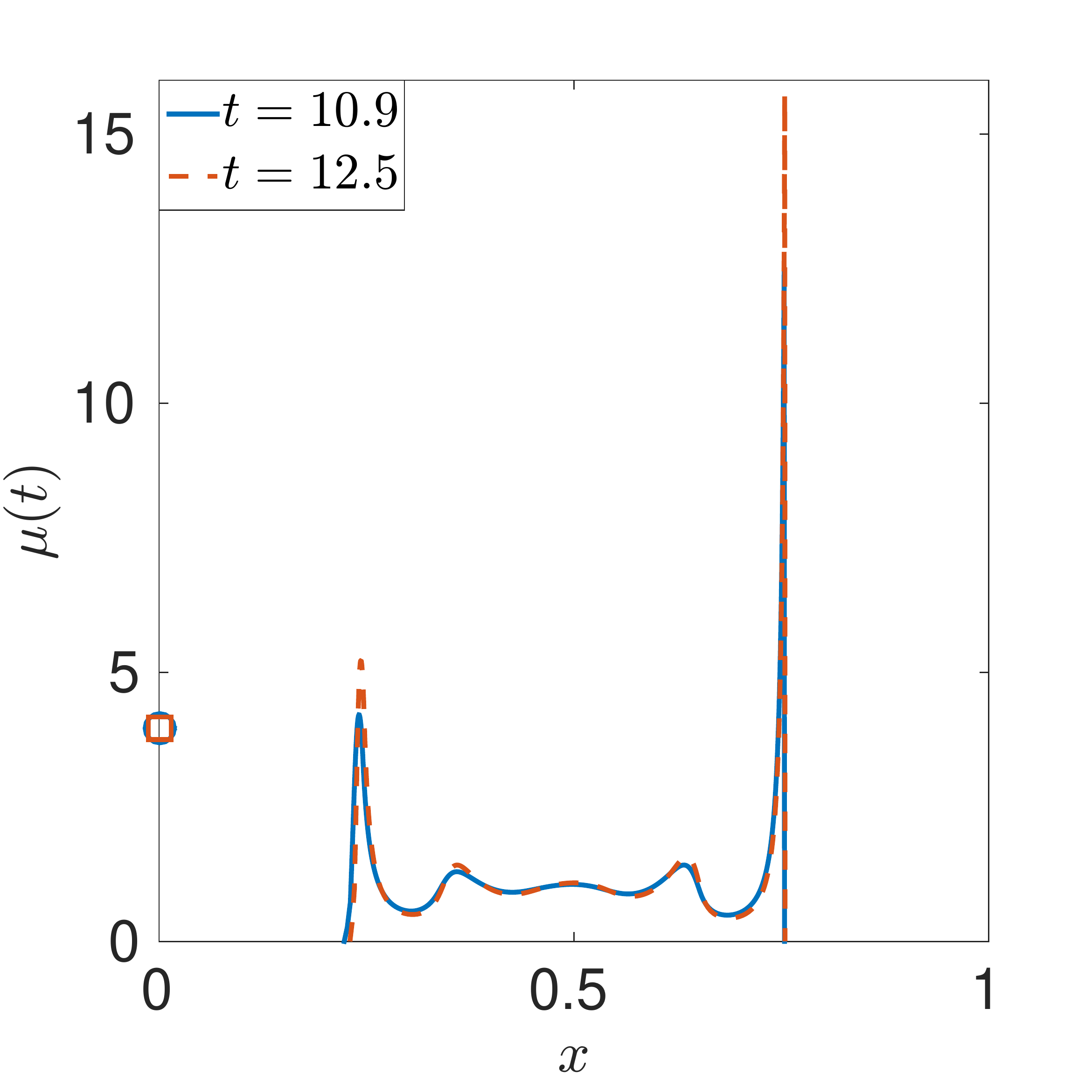} \\
(a) & (b)
 \end{tabular}
\caption{Simulations with potential \eqref{eqn:C2intpot} showing the first mass transfer from the boundary to the free swarm in the diffusive model. (a) Snapshots of the diffusive model \eqref{eqn:eps-eqn} with $\nu = 10^{-7}$. (b) Snapshots of the plain aggregation model \eqref{eqn:agg-eqn2}. Concentrations are represented as circle and square markers for $t = 10.9$ and $t=12.5$ respectively. The masses of concentrations have been magnified $10$ times for clarity.}
\label{fig:SmoothstatesTran}
\end{center}
\end{figure}

Figure~\ref{fig:SmoothstatesTran}(a) shows the onset of the mass transfer at $t=10.9$ for simulations with $\nu = 10^{-7}$. Also, by $t=12.5$ we see that mass has elongated away from the origin and has begun forming a new bump. The transfer of mass occurs repeatedly in the diffusive model as the solution evolves further through time, though less mass is transferred each time. Generally, this mass will either form a new bump or join with the next nearest bump.  In contrast, Figure~\ref{fig:SmoothstatesTran}(b) shows that in the plain aggregation model the concentration at the origin does not change and the five bumps in the free swarm just become sharper, as they tend toward five delta concentrations.

%

Mass transfers are tightly linked to the $2$-Wasserstein distance between solutions $\mu_\nu(t)$ and $\mu(t)$, as well as to the energy evolution of the diffusive model. 
Figure~\ref{fig:SmoothW2FVMvsP}(a) shows the evolution in time of the distance between the two solutions, for various $\nu$. In each plot we see a significant increase in the growth of the distance at exactly the times when mass first transfers away from the origin. Additionally we observe that decreasing $\nu$ keeps the distance between solutions smaller for longer times. 

Complementary to looking at the $2$-Wasserstein distance, Figure~\ref{fig:SmoothW2FVMvsP}(b) compares the energies of the diffusive model and of the plain aggregation model (see \eqref{eqn:reg_energy} and \eqref{eqn:energy}). We observe that the energies are close, again, up until the first mass transfer occurs. 
The energy plots also show that the diffusive model enables solutions to reach lower energies where the plain aggregation model gets stuck at a higher energy that corresponds to an energetically unstable steady state.

\begin{figure}[htb!]
  \begin{center}
   \begin{tabular}{cc}
 \includegraphics[width=0.46\textwidth]{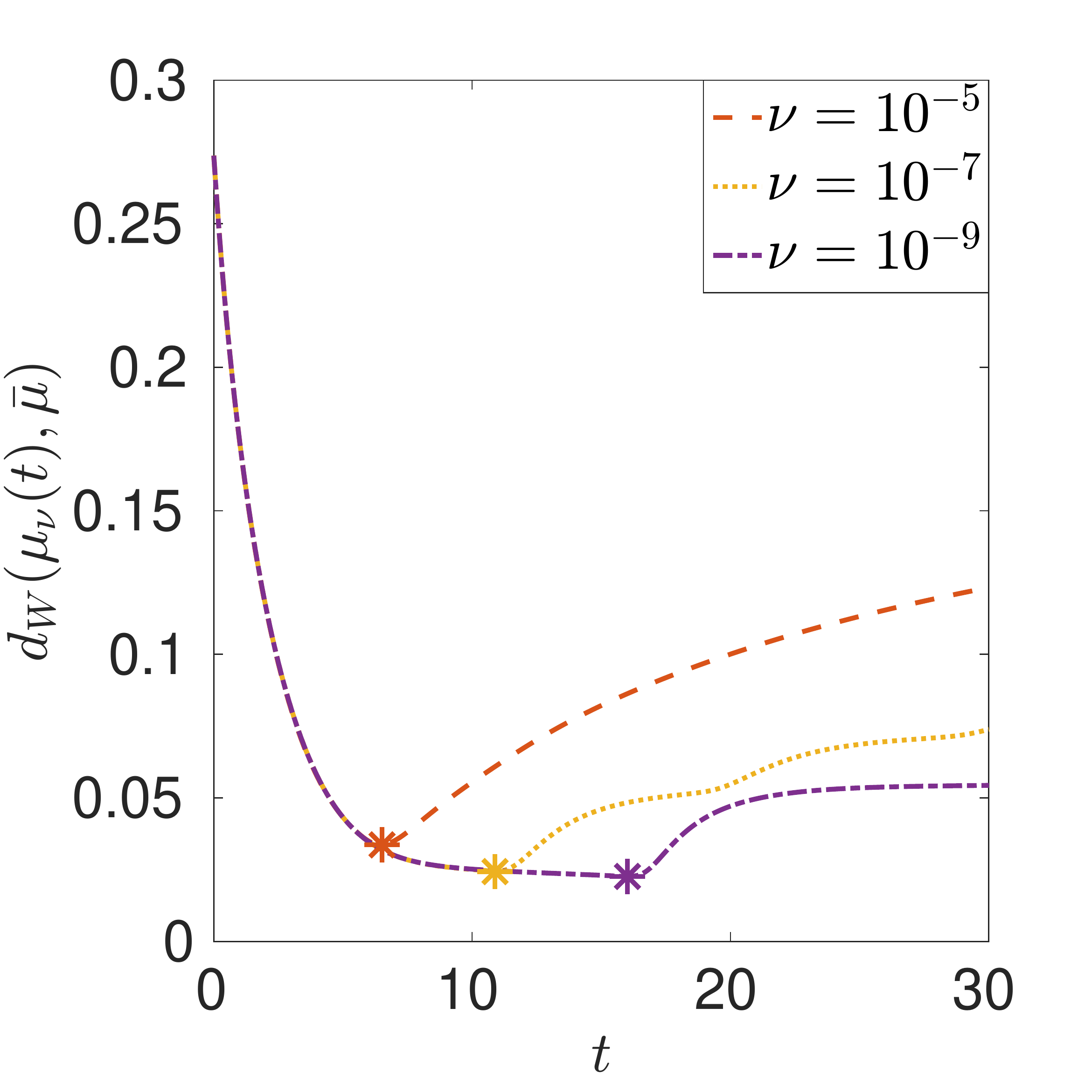} &
 \includegraphics[width=0.46\textwidth]{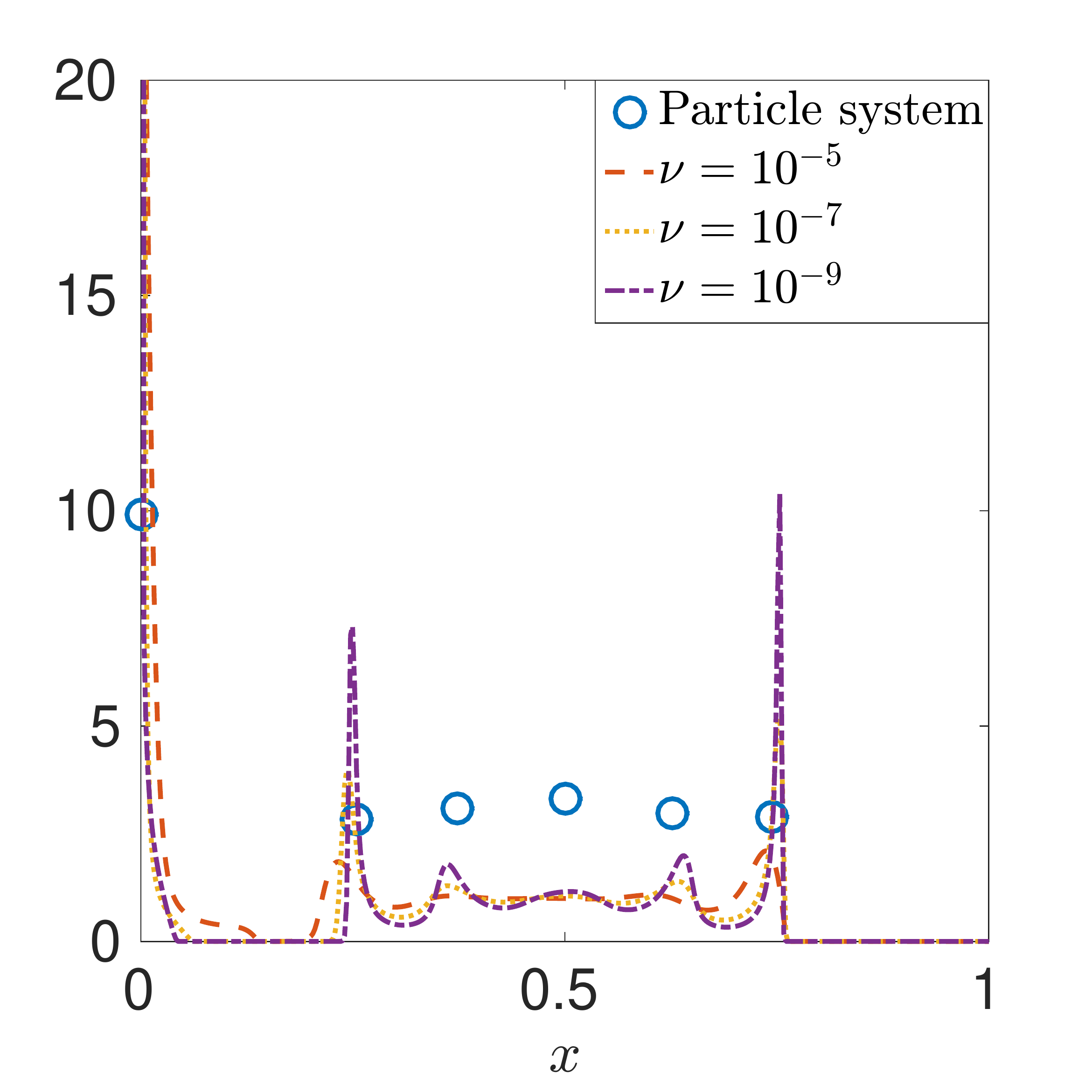} \\
 (a) & (b) 
 \end{tabular}
\caption{Results with potential \eqref{eqn:C2intpot}. (a) $2$-Wasserstein distance between solutions to the diffusive model and the (unstable) equilibrium of the plain aggregation model, for various choices of $\nu$. Markers have been placed at  $t = 6.5$, $t=10.9$ and $t=16$ for $\nu = 10^{-5}$, $\nu = 10^{-7}$ and $\nu = 10^{-9}$ respectively, corresponding to the times of the first mass transfer (see also Figure \ref{fig:SmoothW2FVMvsP}); these times also correspond to when $\mu_\nu(t)$ is closest to $\bar{\mu}$. (b) Solutions to the diffusive model at the times marked in (a). The circles represent concentrations of the equilibrium $\bar{\mu}$, where they have been magnified $25$ times for clarity.}
\label{fig:SmoothW2closestStates}
\end{center}
\end{figure}

With no mechanism to break apart the delta concentration at the boundary (see Figure \ref{fig:SmoothstatesTran}(b)), the plain aggregation model evolves into a steady state $\bar{\mu}$ that consists of six delta concentrations (one at origin and five in the interior) -- see Figure \ref{fig:SmoothW2closestStates}(b). This steady state is {\em not} a minimizer of energy \eqref{eqn:energy}. This flaw of the plain aggregation model was pointed out in \cite{FeKo2017}. On the other hand, the mechanism of mass transfer in the diffusive model enables solutions $\mu_\nu(t)$ to bypass the unstable equilibrium $\bar{\mu}$ of the plain aggregation equation. 

Figure~\ref{fig:SmoothW2closestStates}(a) shows the $2$-Wasserstein distance $d_W(\mu_\nu(t),\bar{\mu})$ between the solutions of the diffusive model and the unstable equilibrium of the plain aggregation model. The plots show that $d_W(\mu_\nu(t),\bar{\mu})$ achieves its minimum (i.e.,  the diffusive model comes nearest to the plain aggregation steady state) exactly at the times of the first mass transfer.  The solutions $\mu_\nu(t)$ at these times, consisting of multiple smoothed aggregates, are shown in Figure~\ref{fig:SmoothW2closestStates}(b). One can see indeed that the smaller the $\nu$ the closer the diffusive solutions pass by the plain aggregation equilibrium. 

It is expected that throughout time evolution, solutions $\mu_\nu(t)$ of the diffusive model bypass other unstable equilibria of the plain aggregation model. This can be observed for instance in the staircase-like evolution of the energy for $\nu=10^{-7}$ in Figure \ref{fig:SmoothW2FVMvsP}(b) (dotted line). The various plateaus of the energy correspond exactly to solutions being temporarily trapped near an unstable equilibria of the plain aggregation model (one could think of these configurations as metastable states for the diffusive model), while the drops in energy correspond to mass transfers.



\smallskip

\begin{figure}[htb!]
  \begin{center}
  \begin{tabular}{cc}
 \includegraphics[width=0.46\textwidth]{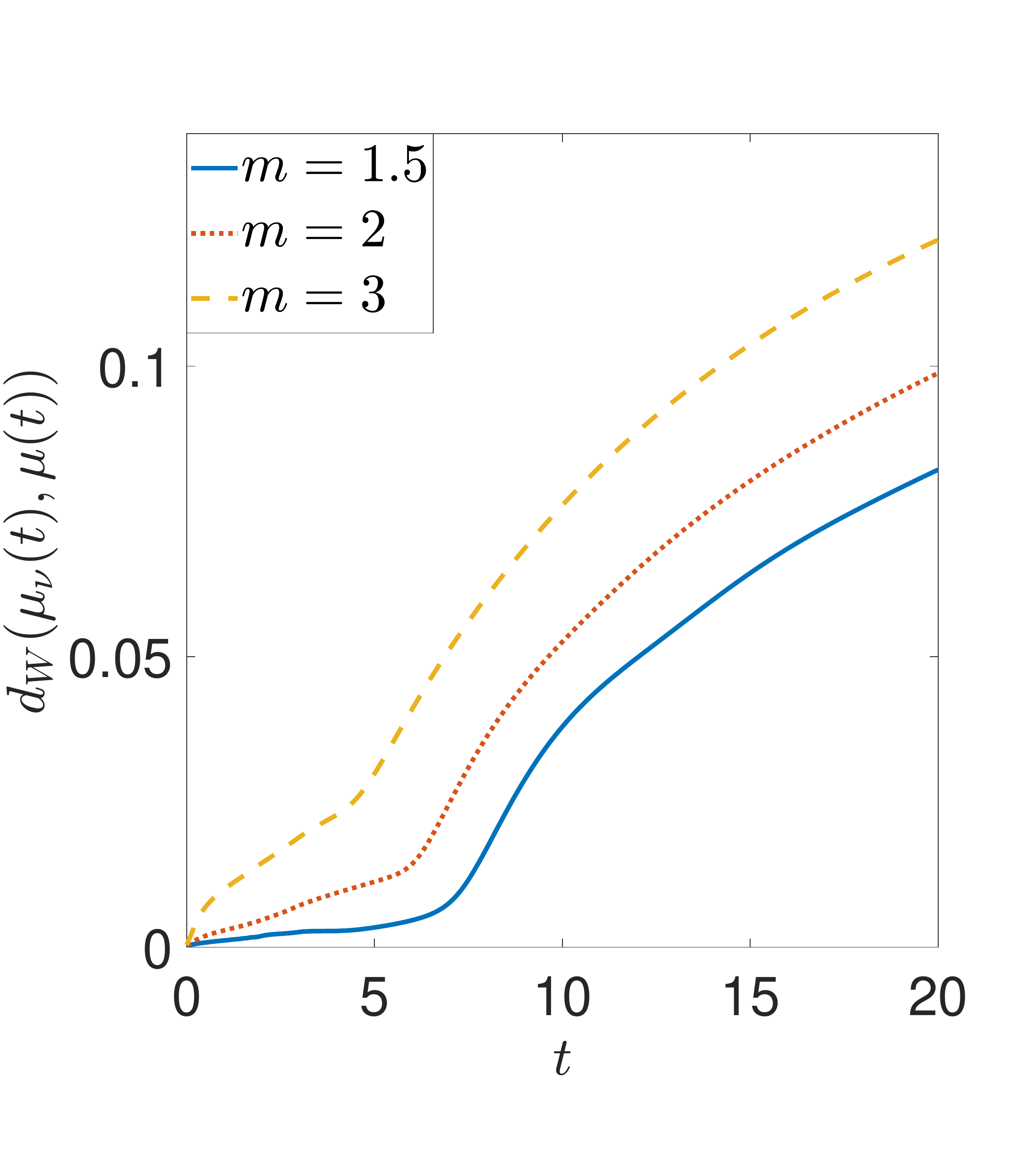} &
 \includegraphics[width=0.46\textwidth]{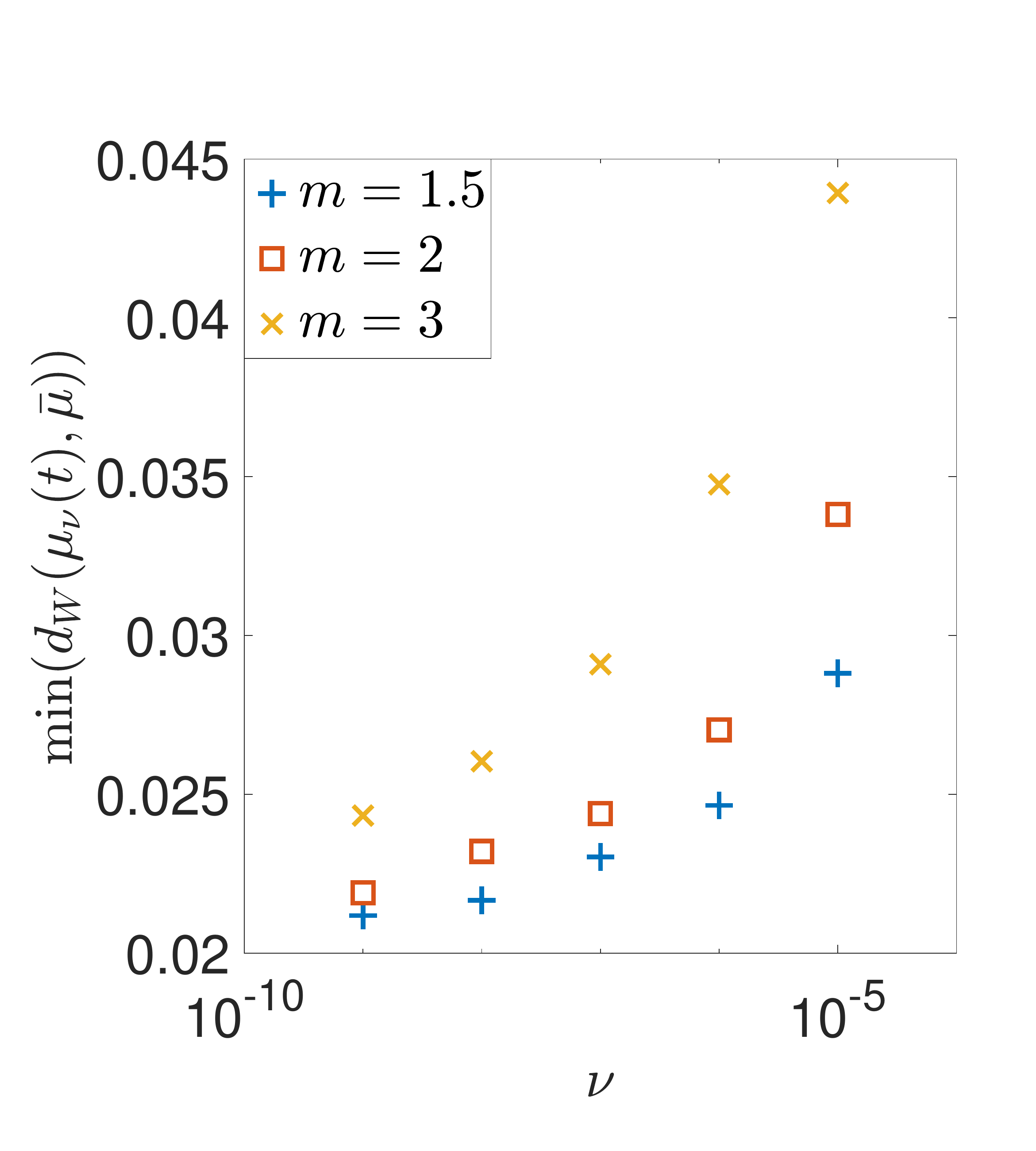} \\
 (a) & (b)
  \end{tabular}
\caption{Results for various values of the exponent $m$. (a) Decreasing $m$ improves the approximation by nonlinear diffusion for fixed $\nu$ (here $\nu=10^{-5}$). (b) By decreasing $m$, the solutions of the diffusive model pass more closely by the unstable equilibrium of the plain aggregation model. }
\label{fig:Varym}
\end{center}
\end{figure}

{\em Simulations with other values of exponent $m$.} The results above have been confirmed with other values of $m$ as well. We performed numerical simulations with $m=1.5$ and $m=3$ starting from the same initial density $\mu^0$. The most notable distinction in this regard is the rate of convergence that depends on $m$.  Specifically, note that in Theorem 1.1 we have derived an explicit upper bound for the convergence rate of $d_W^2(\mu_\nu(t),\mu(t))$ at fixed times, as $\nu \to 0$. The rate (not necessarily sharp) is $\nu^{\beta}$, where $\beta$ depends on $m$ and the dimension $d$ such that 
\[
m_1 < m_2 \quad \implies \quad \beta_1 \geq \beta_2.
\]
Therefore, we expect better convergence at fixed times for lower values of $m$. A numerical validation of this fact is shown in Figure \ref{fig:Varym}(a) for $\nu = 10^{-5}$.

By reproducing the analogues of Figures \ref{fig:SmoothstatesEarly}-\ref{fig:SmoothW2closestStates} for $m=1.5$ and $m=3$ we found indeed that the approximation by nonlinear diffusion gets better/sharper by decreasing $m$. For instance: in Figure \ref{fig:SmoothstatesEarly}(a) the boundary layer at origin gets steeper and narrower with decreasing $m$ and hence it better approximates the Dirac accumulation at the origin in the plain aggregation model.  Aside from the general observation made above (rate of convergence $\nu^{\beta}$ improves with lowering $m$), a formal argument for this fact is that at high concentrations $\rho$, the diffusion $\rho^m$ decreases with $m$. Also, at later times  diffusive solutions with lower values of $m$ capture more sharply the interior delta aggregations of the plain aggregation model (case $m=2$ is shown in Figure \ref{fig:SmoothW2closestStates}(b)). 

Finally, to conclude the discussion on varying $m$,  the first mass transfer from the boundary occurs faster for larger values of $m$ (where there is more diffusion).  In energy plots such as that shown in Figure \ref{fig:SmoothW2FVMvsP}(b), the energy staircasing gets accelerated by increasing $m$ (more diffusion, faster mass transfers). Also,  the minimum distance between the solutions of the diffusive model and the equilibrium of the plain aggregation model occurs faster for larger $m$. On the other hand, these minimum distances decrease with lowering $m$ -- see Figure \ref{fig:Varym}(b). A detailed account on the simulations with various $m$ will be included in the upcoming PhD thesis of one the authors \cite{KovacicPhD2018}. 


\subsection{Time evolution: $C^0$ interaction potential}\label{subsec:c0intpot}

We consider now the $C^0$ interaction potential \eqref{eqn:C0intpot} instead of its regularized, smooth version \eqref{eqn:C2intpot}. While potential \eqref{eqn:C0intpot} does not satisfy the conditions in Theorem~\ref{thm:conv}, we show that some of the same general observations made in Section~\ref{subsec:c2intpot} can be made in this case as well.

\begin{figure}[htb!]
  \begin{center}
  \begin{tabular}{cc}
 \includegraphics[width=0.46\textwidth]{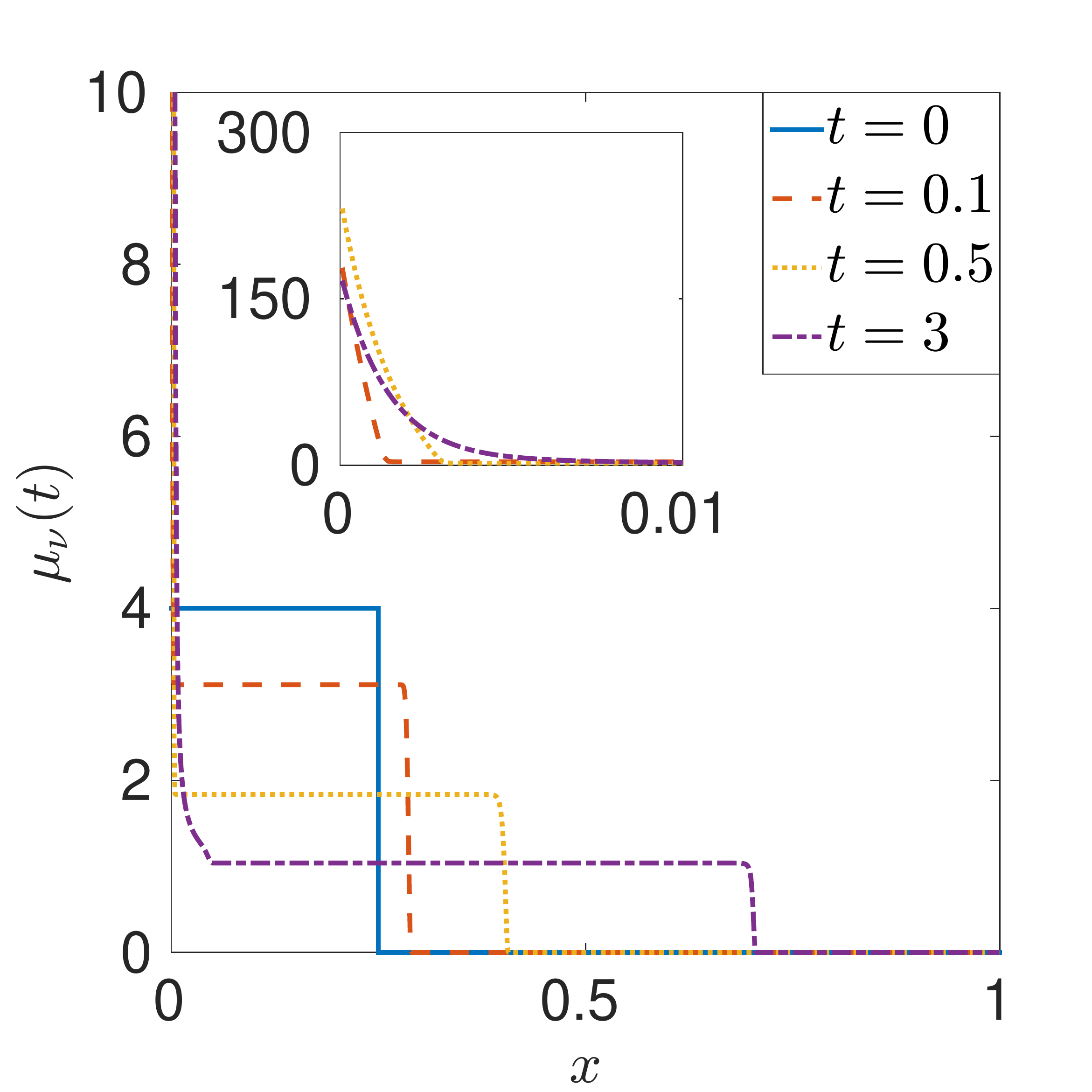} &
 \includegraphics[width=0.46\textwidth]{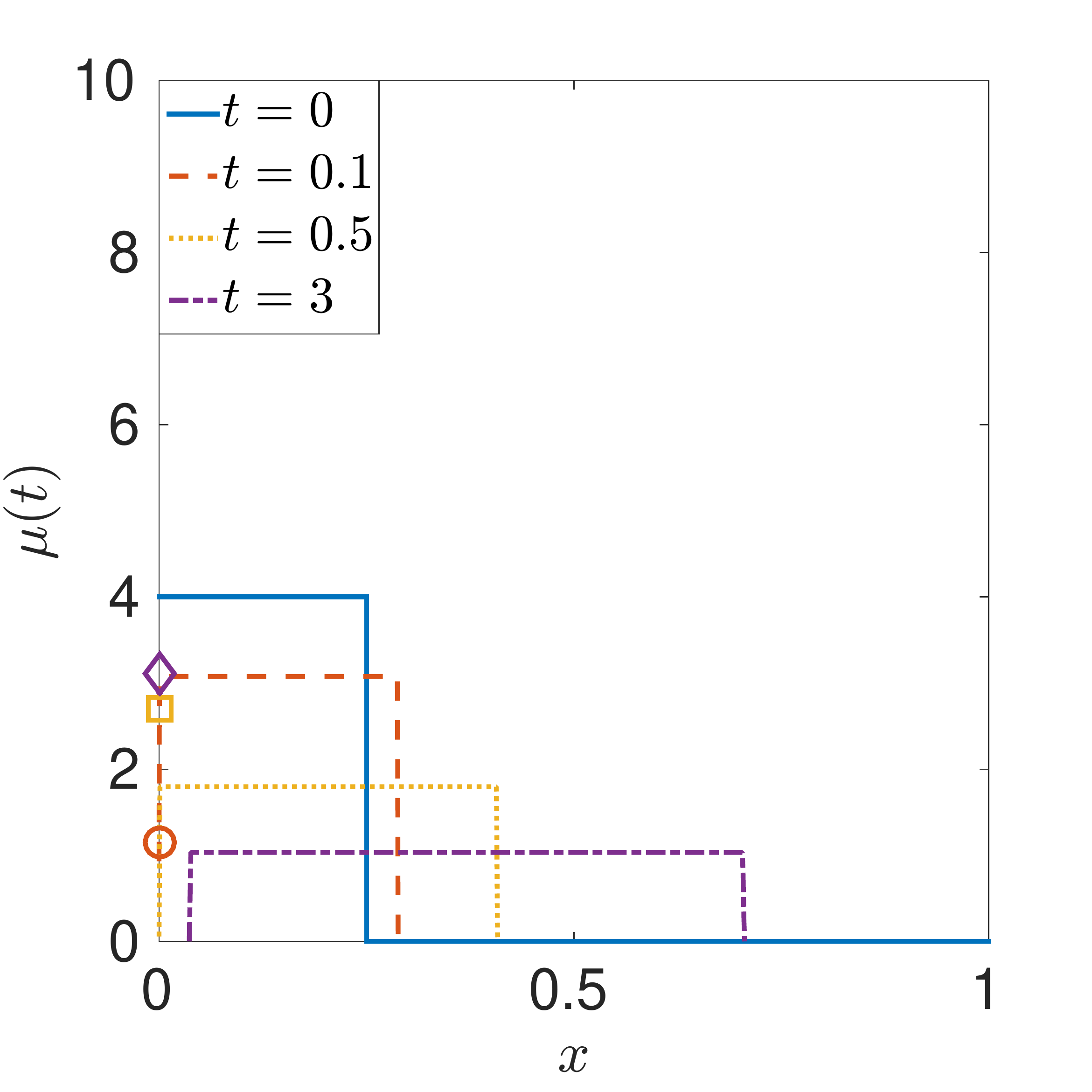} \\
 (a) & (b)
  \end{tabular}
\caption{Simulations with potential \eqref{eqn:C0intpot} showing early time dynamics. (a) Snapshots of the diffusive model \eqref{eqn:eps-eqn} with $\nu = 10^{-6}$. An insert has been included to show the layer of mass near the origin more clearly. (b) Snapshots of the plain aggregation model \eqref{eqn:agg-eqn2}. Concentrations are represented as circle, square, and diamond markers for $t = 0.1$, $t = 0.5$, and $t=3$ respectively. The masses of concentrations have been magnified 10 times for clarity.}
\label{fig:NonsmoothstatesEarly}
\end{center}
\end{figure}

For early times we find that the solutions of the diffusive and plain aggregation models remain qualitatively similar (see Figure~\ref{fig:NonsmoothstatesEarly}): the initial density moves apart, with some mass staying near, or on, the origin and the rest spreading away from the wall. As in Section~\ref{subsec:c2intpot} we find again that the diffusive model has a thin, sharp layer of mass near the origin where the plain aggregation model concentrates mass exactly at the origin. One notable difference from Section~\ref{subsec:c2intpot} is that now the diffusive solution $\mu_\nu$ consists of only a single component for all tests we performed.

\begin{figure}[thb!]
  \begin{center}
  \begin{tabular}{cc}
 \includegraphics[width=0.46\textwidth]{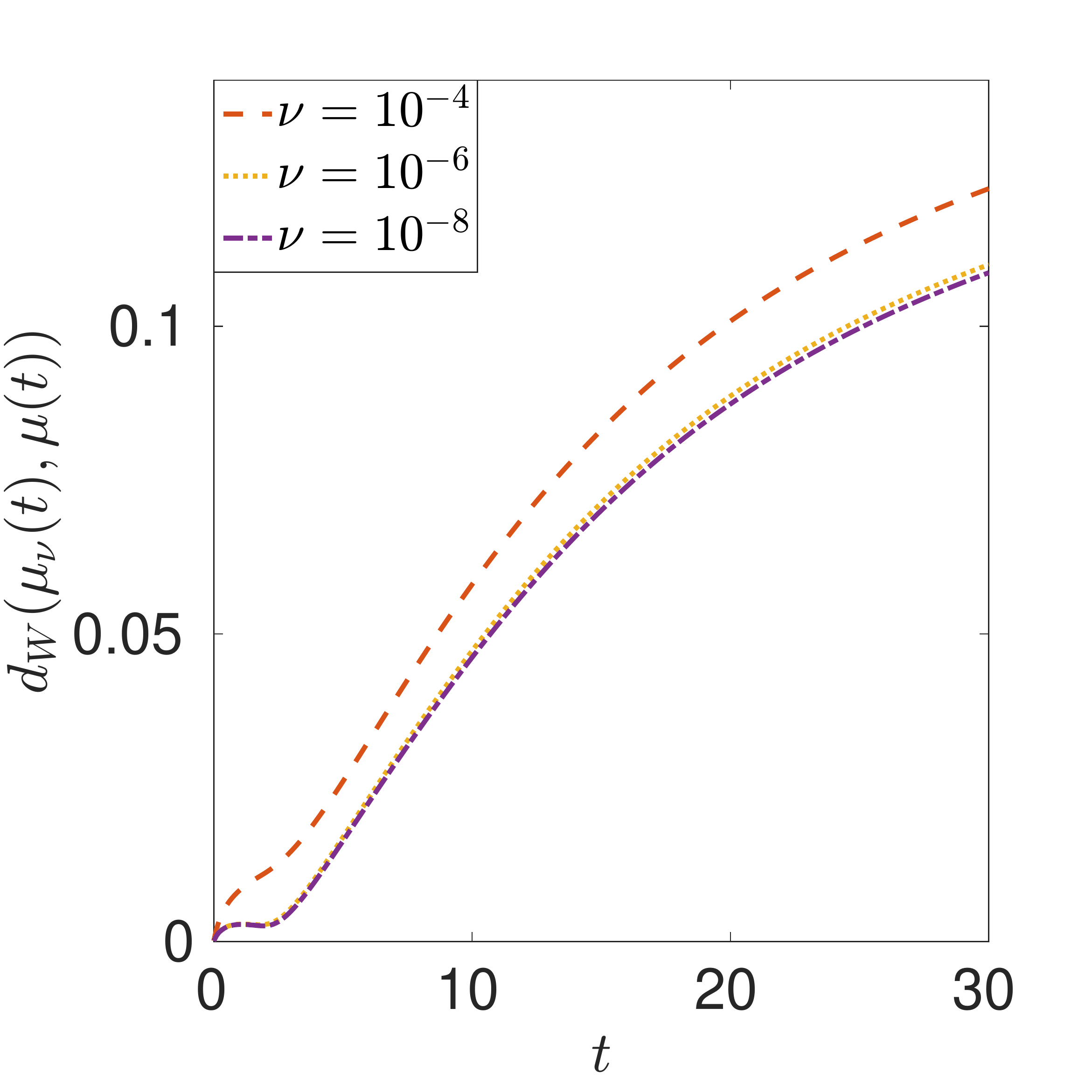} &
 \includegraphics[width=0.46\textwidth]{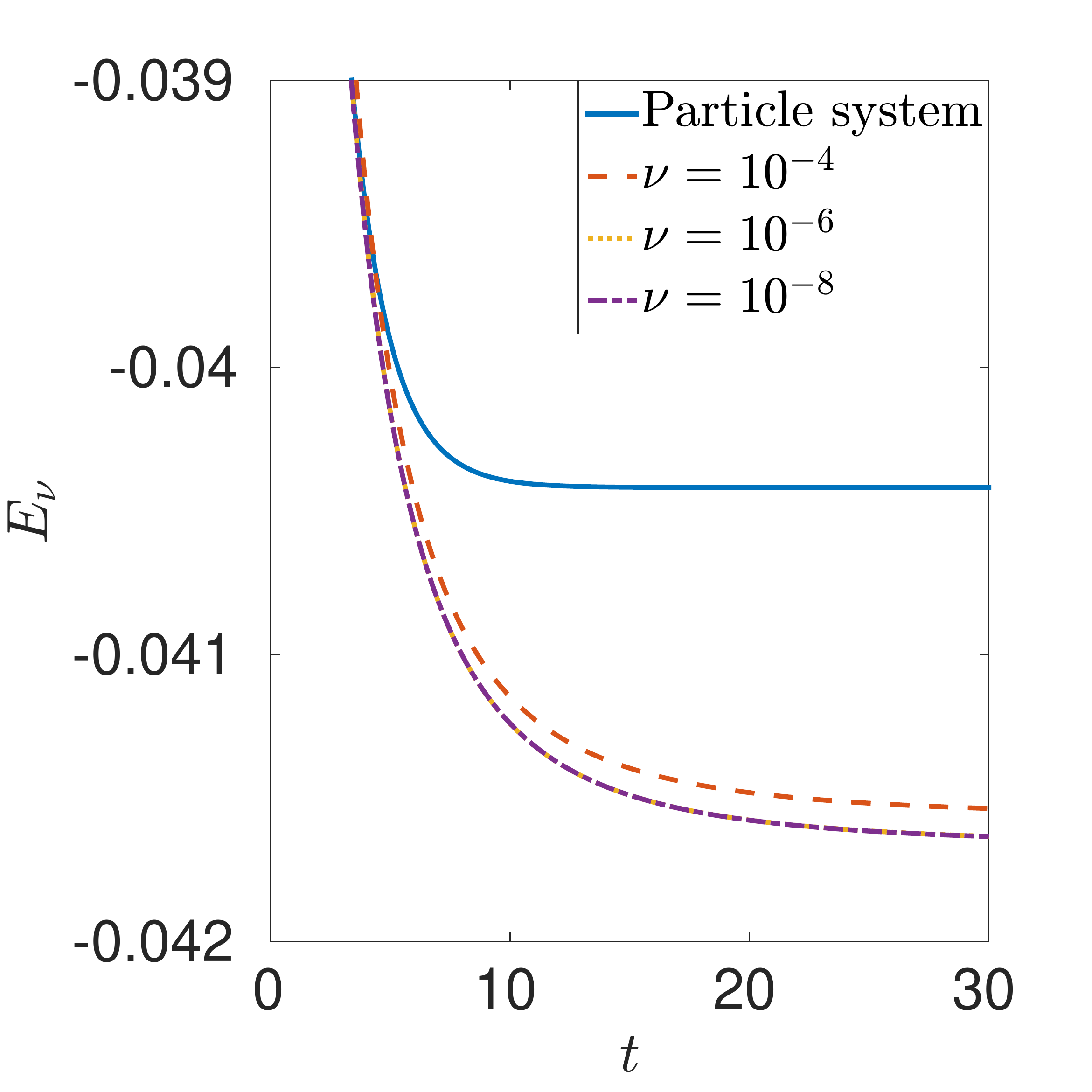} \\
 (a) & (b) 
 \end{tabular}
\caption{Results with potential \eqref{eqn:C0intpot}. (a) $2$-Wasserstein distance between the diffusive and plain aggregation solutions for various choices of $\nu$. (b) Energy \eqref{eqn:reg_energy} of solutions to the diffusive model through time for various choices of $\nu$. Also included is the energy \eqref{eqn:energy} of the solution to the particle model through time (solid line).}
\label{fig:NonsmoothW2FVMvsP}
\end{center}
\end{figure}

We also consider the $2$-Wasserstein distance $d_W(\mu_\nu(t),\mu(t))$ between the two solutions and find similar results to  Section~\ref{subsec:c2intpot} -- see Table~\ref{tab:NonsmoothW2dist} and Figure \ref{fig:NonsmoothW2FVMvsP}(a). Specifically, we find that for fixed times the distance decreases as $\nu$ is decreased toward zero, suggesting that Theorem \ref{thm:conv} may be refined to include less regular interaction potentials. 
We also find that the distance between the two solutions increases as time goes forward, with a faster growth rate at later times.
\begin{table}
\begin{centering}
\begin{tabular}{|c|c|c|c|}
\hline 
$\nu$ & $t=0.1$ & $t=0.5$ & $t=3$ \\ 
\hline 
$10^{-3}$ & $6.8548\scie{-3}$ & $2.4142\scie{-2}$ & $6.6831\scie{-2}$ \\ 
\hline 
$10^{-4}$ & $2.9493\scie{-3}$ & $1.1424\scie{-2}$ & $4.3318\scie{-2}$ \\ 
\hline 
$10^{-5}$ & $2.3620\scie{-3}$ & $8.9352\scie{-3}$ & $3.7865\scie{-2}$ \\ 
\hline 
$10^{-6}$ & $2.3166\scie{-3}$ & $8.6382\scie{-3}$ & $3.6588\scie{-2}$ \\ 
\hline 
$10^{-7}$ & $2.3161\scie{-3}$ & $8.6054\scie{-3}$ & $3.6235\scie{-2}$ \\ 
\hline 
\end{tabular}
\captionof{table}{$2$-Wasserstein distance $d_W(\mu_\nu(t),\mu(t))$ between solutions of the diffusive model and solutions of the plain aggregation model for various choices of $\nu$ at some early times.} \label{tab:NonsmoothW2dist} \par
\end{centering}
\end{table}
\medskip

Figure~\ref{fig:NonsmoothW2FVMvsP}(b) shows the energies of the diffusive and plain aggregation solutions. We find again that with diffusion, solutions achieve states of lower energy than the plain aggregation model. We do not see the same energy staircase pattern as in Figure~\ref{fig:SmoothW2FVMvsP}(b) however. This is not unexpected actually, as the reason we observed the staircase pattern in Section~\ref{subsec:c2intpot} was because the diffusive solution consisted of multiple disjoint components. The staircasing was highly linked to instances of mass from the origin gradually pulling away, leaving the origin, and moving to join the free swarm. Since with potential \eqref{eqn:C0intpot} the diffusive solution does not form multiple disjoint components, there is no mass transfer and hence the mechanism for energy staircasing is missing. 

\begin{figure}[thb!]
  \begin{center}
  \begin{tabular}{cc}
 \includegraphics[width=0.46\textwidth]{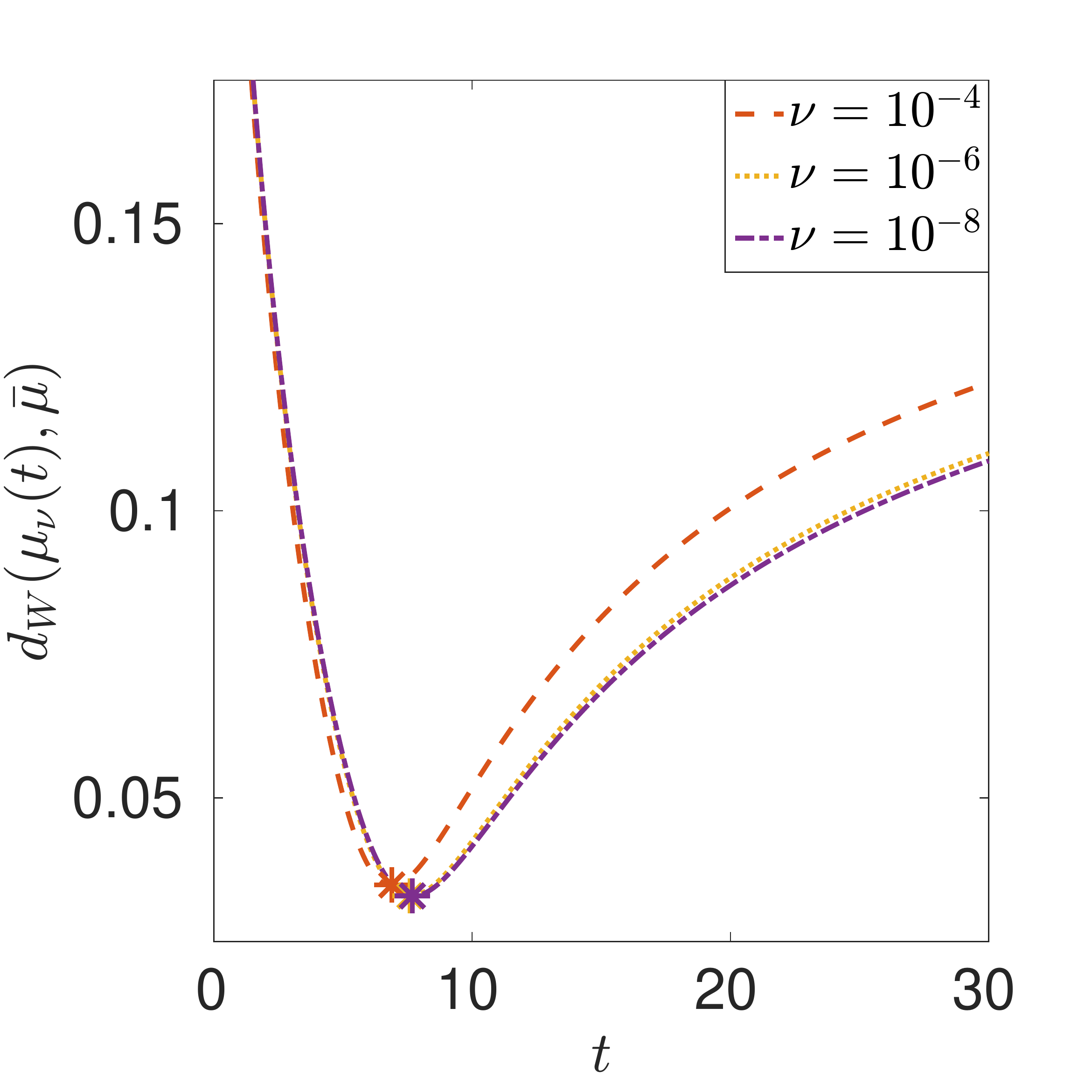} &
 \includegraphics[width=0.46\textwidth]{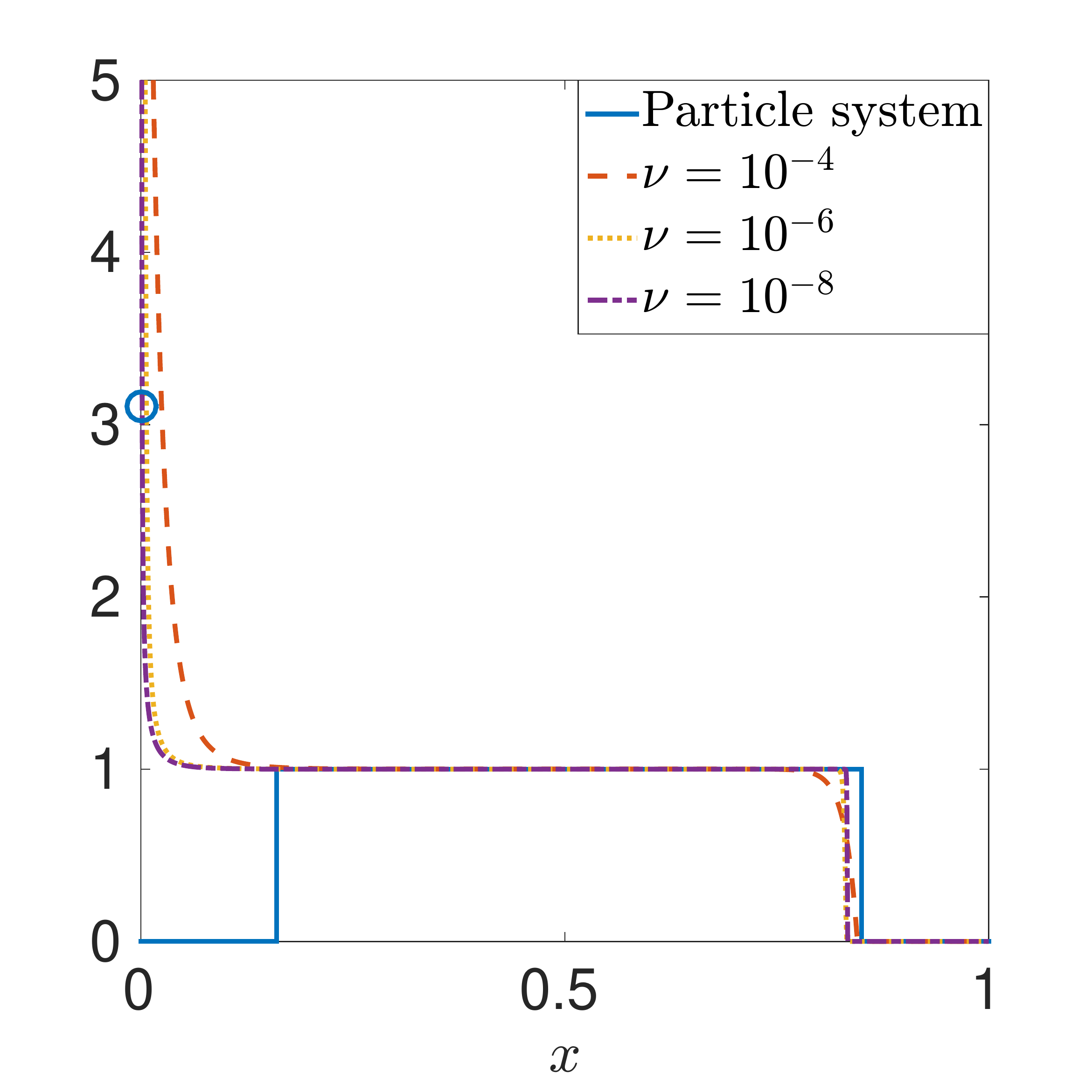} \\
 (a) & (b)
 \end{tabular}
\caption{Results with potential \eqref{eqn:C0intpot}. (a) $2$-Wasserstein distance between solutions to the diffusive model and the (unstable) equilibrium of the plain aggregation model, for various choices of $\nu$. Markers have been placed at $t = 6.9$, $t=7.6$, and $t=7.9$ for $\nu = 10^{-4}$, $\nu = 10^{-6}$, and $\nu = 10^{-8}$ respectively, corresponding to the times when $d_W(\mu_\nu(t),\bar{\mu})$ achieves its minimum. (b) Solutions to the diffusive model at the times marked in (a), respectively for each $\nu$. The solid line and the circle marker at origin (indicating a delta concentration) represents the unstable equilibrium $\bar{\mu}$ of the plain aggregation model. The concentration has been magnified $10$ times for clarity.}
\label{fig:NonsmoothW2closestStates}
\end{center}
\end{figure}

We also find again that the diffusive model bypasses the unstable equilibrium of the plain aggregation model. Figure~\ref{fig:NonsmoothW2closestStates}(a) shows that the solutions $\mu_\nu(t)$ of the diffusive model come near to the unstable equilibrium $\bar{\mu}$ of the plain aggregation model and decreasing $\nu$ causes this distance to decrease, though not as noticeably as for the smooth potential (Figure~\ref{fig:SmoothW2closestStates}(a)). We have also included markers at times where $d_W(\mu_\nu(t),\bar{\mu})$ achieves its minimum though we do not see these times being significant to $d_W(\mu_\nu(t),\mu(t))$ in Figure~\ref{fig:NonsmoothW2FVMvsP}(a) or to the energies of the diffusive model in Figure~\ref{fig:NonsmoothW2FVMvsP}(b). We believe that, as for the energy staircasing, it is the lack of mass transfer of solutions that is the cause here. We also think that for the same reason the curves in Figures~\ref{fig:NonsmoothW2FVMvsP}(a) and \ref{fig:NonsmoothW2closestStates}(a) are not as differentiated as for the smooth potential \eqref{eqn:C2intpot}.

Figure~\ref{fig:NonsmoothW2closestStates}(b) shows the solutions $\mu_\nu(t)$ compared to the (unstable) equilibrium $\bar{\mu}$ of the plain aggregation model at times corresponding to the minima of $d_W(\mu_\nu(t),\bar{\mu})$, specifically $t = 6.9$, $t=7.6$, and $t=7.9$ for $\nu = 10^{-4}$, $\nu = 10^{-6}$, and $\nu = 10^{-8}$ respectively. Note that the solutions $\mu_\nu(t)$ match up with most of the free swarm component of $\bar{\mu}$. The major qualitative difference is again, that solutions of the diffusive model consist of a single component where the plain aggregation equilibria is formed of two disjoint parts. 

We also note that we performed numerical simulations with potential \eqref{eqn:C0intpot} for other values of the exponent $m$, and found that similar considerations as above hold as well. More details on such simulations can be found in \cite{KovacicPhD2018}.

\subsection{Convergence of energy minimizers}

We now present some numerical evidence for Theorem~\ref{thm:conv_min} with quadratic diffusion ($m=2$) and interaction potential \eqref{eqn:C0intpot}. Before we begin it is important to point out that while the singular potential \eqref{eqn:C0intpot} does not satisfy the conditions of Theorem~\ref{thm:conv}, it does satisfy the conditions of Theorem~\ref{thm:conv_min}. The benefits of using quadratic diffusion and this interaction potential is that we can explicitly calculate the equilibria of the diffusive model and then take their zero diffusion limit.  To find equilibria and local minimizers we will be relying on the variational framework set up in \cite{BeTo2011}, adapted to interaction energies that also contain diffusive terms.

In \cite{FeKo2017} the authors studied equilibria of the plain aggregation model on the half-line $[0,\infty)$ with interaction potential \eqref{eqn:C0intpot}. It was found that in the absence of an external potential, the minimizing equilibrium $\bar{\mu}^*$ of energy $\E$ (given by \eqref{eqn:energy} with $V$=0) is the same as that for the problem in free space \cite{FeHuKo11}, i.e., 
\begin{equation}\label{eqn:MinPA}
\bar{\mu}^* = \One_{[0,1]}.
\end{equation}
Note that this minimizer is unique up to translation. In addition, there exists a one-parameter family of equilibria $\bar{\mu}$ that are not energy minimizers, composed of a concentration at the origin and a constant-valued, compact component away from the boundary (see solid line and circle marker in Figure~\ref{fig:NonsmoothW2closestStates}(b)). For this reason we added an asterisk superscript to the minimizer in \eqref{eqn:MinPA}, to distinguish it from the other (unstable) equilibria.

We now look for explicit equilibria of the diffusive model with quadratic diffusion ($m=2$) and diffusion exponent $\alpha=1$, with the interaction potential \eqref{eqn:C0intpot} and no external potential. By the results in \cite{BuFeHu14} and \cite{FeRa10}, as well as observations of our numerics in Section~\ref{subsec:c0intpot}, we assume that the equilibrium $\munubar$ is continuous, smooth on its support, and composed of a single component with compact support such that $\supp(\munubar) = [0,L]$.

Following \cite{BeTo2011}, we find equilibria $\munubar$ by looking for critical points of the energy $\E_\nu$. By an immediate calculation (see \cite{BeTo2011}, also \cite{FeKo2017}), one finds from \eqref{eqn:energy} that the first variation of $\E_\nu$ vanishes at $\munubar$, provided
\begin{equation}
\label{eqn:Lnu}
\Lambda_\nu(x) = \lambda_\nu \qquad \text{ for } x \in \supp(\munubar),
\end{equation}
for some $\lambda_\nu \in \R$, where
\begin{equation}\label{eqn:MinLam}
\Lambda_\nu(x) = 2\nu\munubar(x) + K\ast \munubar.
\end{equation}
In brief, $\Lambda_\nu(x)$ can be regarded as  the energy per unit mass felt by a test mass at position $x$ (see \cite{BeTo2011}), while $\lambda_\nu$ arises as a Lagrange multiplier.

Furthermore, if 
\begin{equation}
\label{eqn:Lnu-min}
\Lambda_\nu(x) \geq \lambda_\nu \qquad \text{ for } x \not\in \supp(\munubar),
\end{equation}
then $\munubar$ is a local minimizer. The interpretation of \eqref{eqn:Lnu-min} is that transporting mass from the support of $\munubar$ into its complement increases the total energy \cite{BeTo2011}.

The approach we take is to solve $\Lambda_\nu'(x) = 0$ and $\Lambda_\nu''(x) = 0$ for $x \in \supp(\munubar)$ and then show that these equilibria must necessarily be local minimizers (i.e., satisfy \eqref{eqn:Lnu-min}). We also note that we look for minimizers of unit mass ($\int_0^L \munubar(x) \, dx = 1$) that are continuous as the end of their support:
\begin{equation}\label{eqn:contCond}
\lim_{x \nearrow L} \munubar(x) = 0.
\end{equation}

From \eqref{eqn:MinLam} and \eqref{eqn:C0intpot}, one gets
\begin{equation}\label{eqn:MinLamd2}
\Lambda_\nu''(x) = 2\nu \munubar''(x) + 1- \munubar(x).
\end{equation}
Solving $\Lambda_\nu''(x) = 0$ is trivial and we get the general form of equilibria $\munubar$ to be
\begin{equation}\label{eqn:SSforms}
\munubar(x) = c_1e^{\frac{x}{\sqrt{2\nu}}} + c_2e^{-\frac{x}{\sqrt{2\nu}}} + 1, \qquad x \in [0,L].
\end{equation}
We cannot, however, explicitly solve for the unknowns $c_1$, $c_2$, and $L$ and we resort to obtaining them by numerically solving the remaining three conditions: $\Lambda_\nu' = 0$ in $[0,L]$, the unit mass condition, and continuity at $L$ given by \eqref{eqn:contCond}. We point out that $\Lambda_\nu(x)$ is continuous and that $\Lambda_\nu''(x) = 1 >0$ for $x \not\in \supp(\munubar)$. The strict convexity of $\Lambda_\nu$ outside $\supp(\munubar)$, combined with \eqref{eqn:Lnu}, implies that \eqref{eqn:Lnu-min} necessarily holds, so the equilibria  \eqref{eqn:SSforms} are in fact minimizers of the energy.

\begin{figure}[thb!]
  \begin{center}
  \begin{tabular}{cc}
 \includegraphics[width=0.46\textwidth]{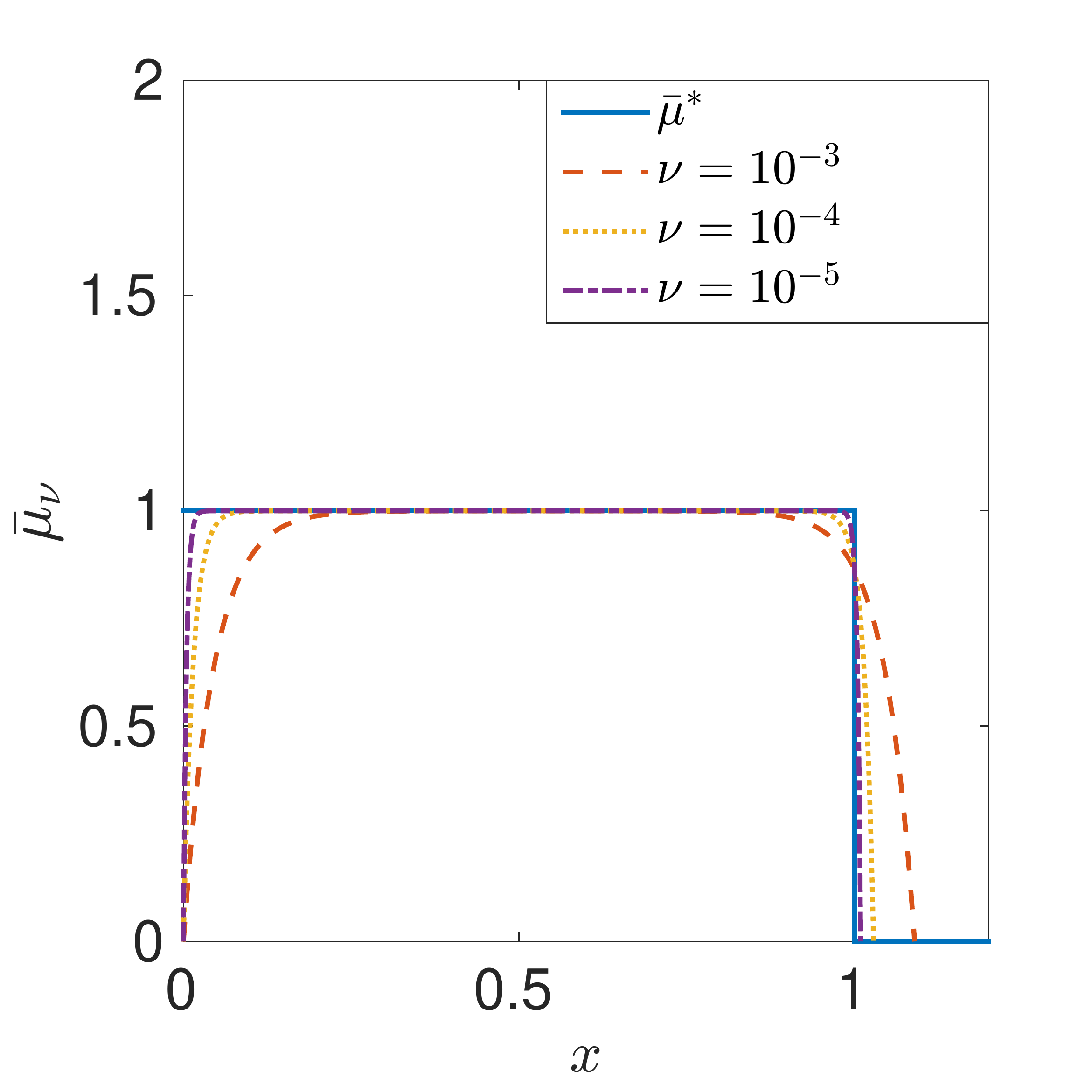} &
 \includegraphics[width=0.46\textwidth]{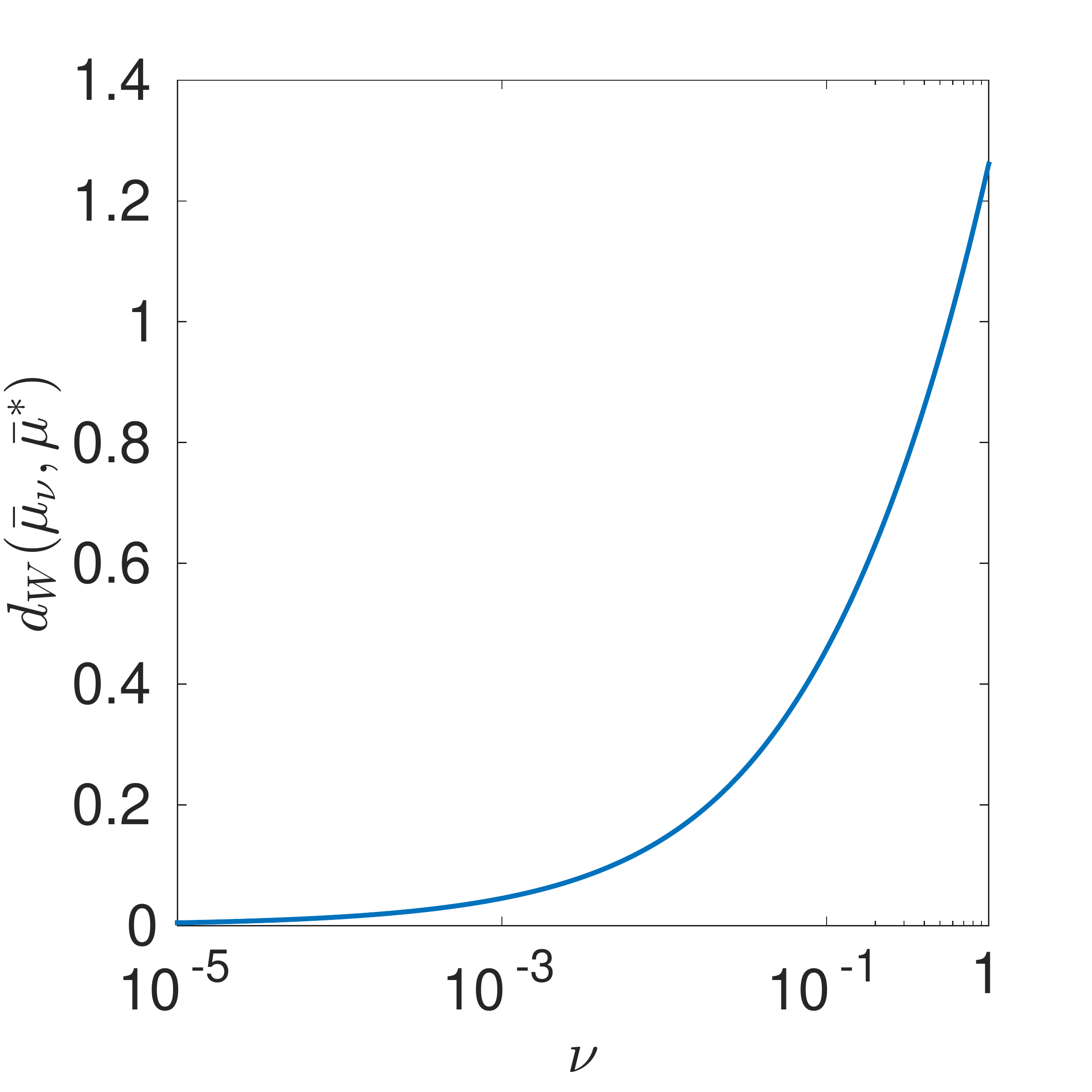} \\
 (a) & (b)
 \end{tabular}
 \end{center}
\caption{(a) Comparison between the energy minimizer $\bar{\mu}^*$  of the plain aggregation model (see \eqref{eqn:MinPA}) and minimizers $\munubar$  of the diffusive model  (see \eqref{eqn:SSforms}) for various $\nu$. (b) The $2$-Wasserstein distance between the minimizers $\munubar$ and $\bar{\mu}^*$  as a function of $\nu$.}
\label{fig:ExpCompare}
\end{figure}

Next we compare the minimizers \eqref{eqn:SSforms} of $\E_\nu$ with the minimizer \eqref{eqn:MinPA} of $\E$ for the plain aggregation model. Figure~\ref{fig:ExpCompare}(a) shows that $\munubar$ qualitatively approach $\bar{\mu}^*$, while Figure~\ref{fig:ExpCompare}(b) provides direct quantitative numerical evidence for Theorem \ref{thm:conv_min}, namely that minimizers of $\E_\nu$ approach (in the $2$-Wasserstein metric) minimizers of $\E$ in the zero diffusion limit.

\begin{figure}[thb!]
  \begin{center}
  \begin{tabular}{cc}
 \includegraphics[width=0.46\textwidth]{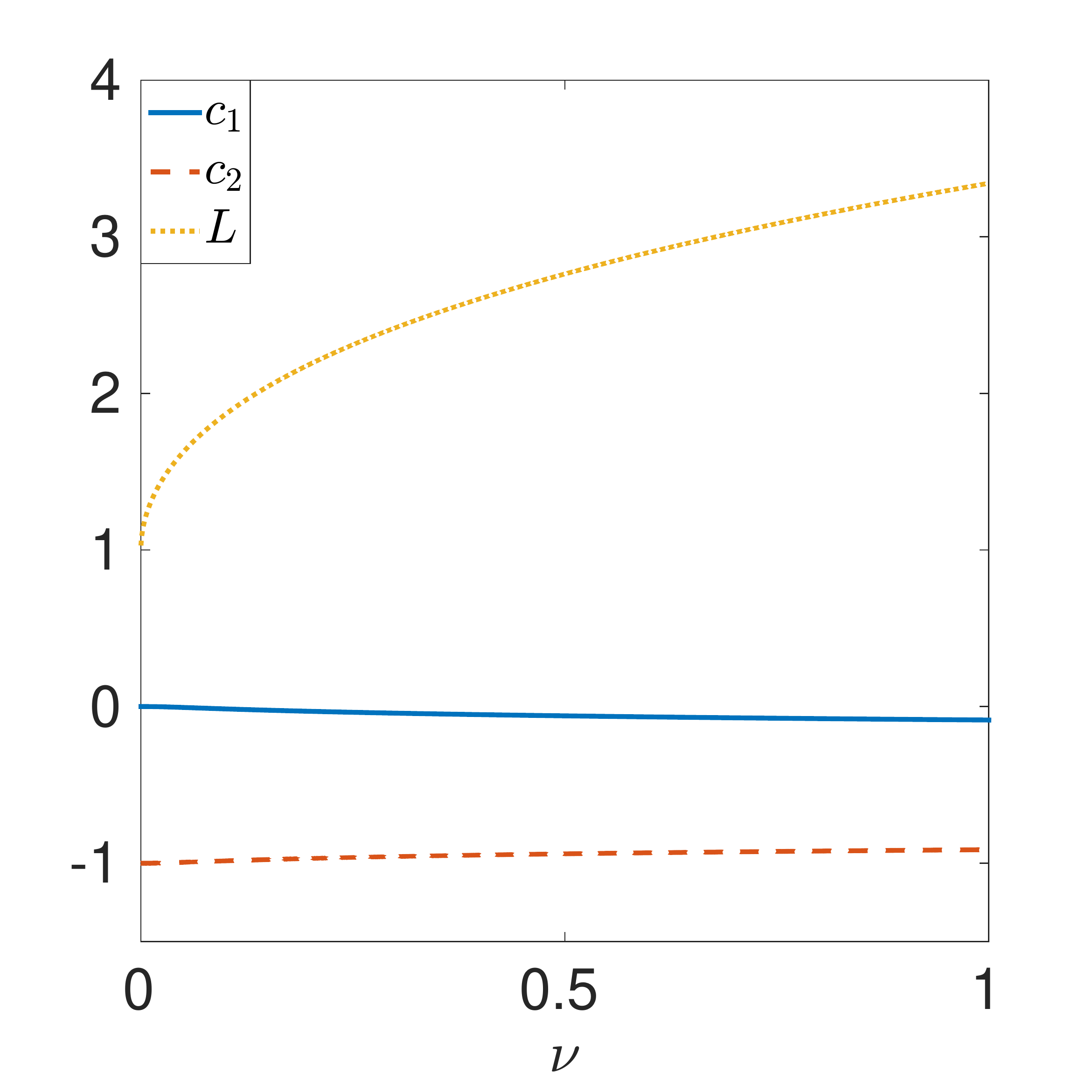} &
 \includegraphics[width=0.46\textwidth]{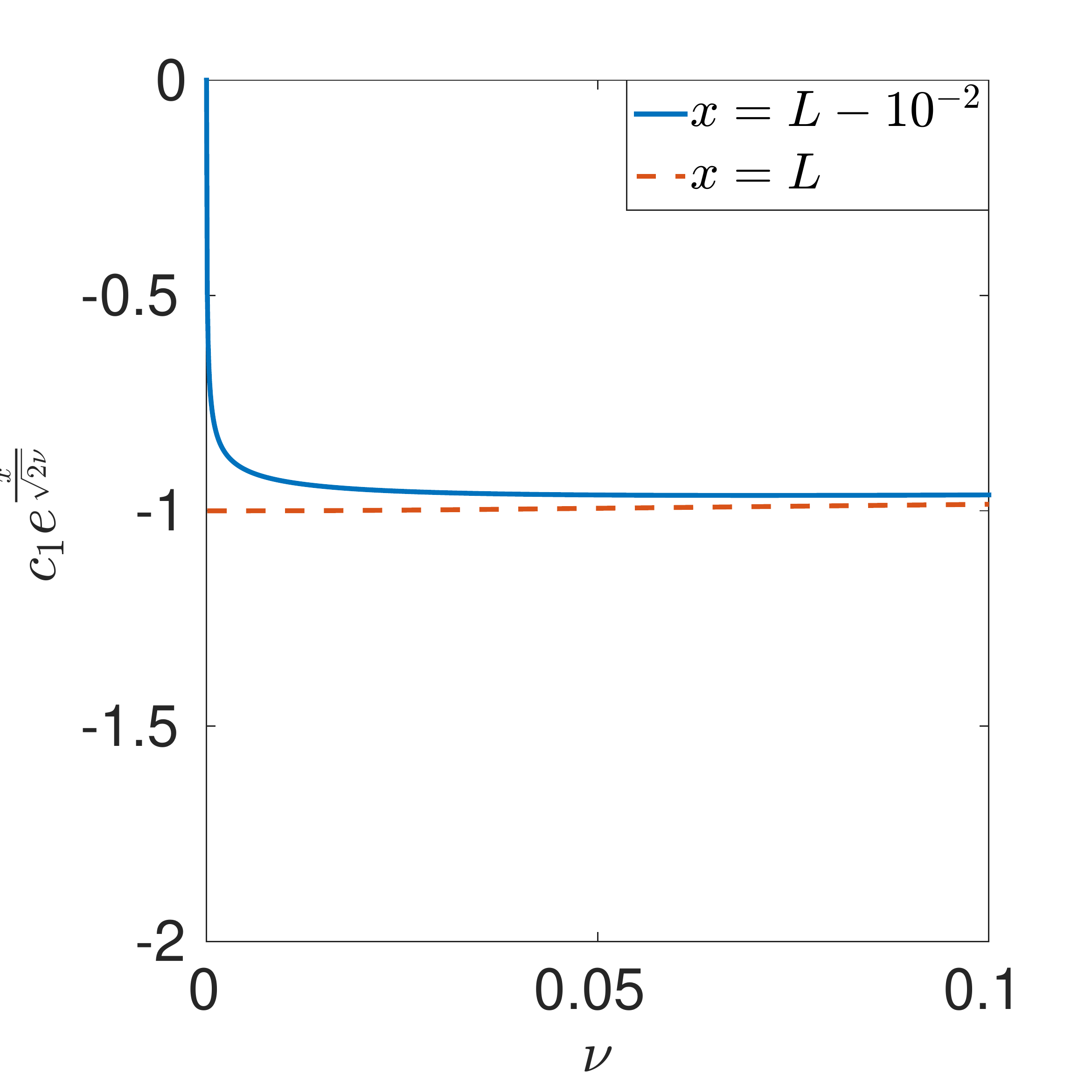} \\
 (a) & (b)
 \end{tabular}
 \end{center}
\caption{(a) Numerically calculated $c_1$, $c_2$, and $L$ as functions of $\nu$. (b) Numerically calculated behaviour of $c_1 e^{\frac{x}{\sqrt{2\nu}}}$ as $\nu \to 0$ for $x$ at $L$ and near $L$.}
\label{fig:MinConsts}
\end{figure}

Figure~\ref{fig:MinConsts}(a) shows the numerically calculated values of $c_1$, $c_2$, and $L$ as functions of $\nu$ (the lowest $\nu$ for which such results have been obtained is $\nu = 10^{-5}$). Observe that $L$ tends to $1$ as $\nu$ tends to zero which is in agreement with the minimizer \eqref{eqn:MinPA} of the plain aggregation model. Furthermore notice that $c_1$ and $c_2$ approach $0$ and $-1$, respectively, in the zero diffusion limit. As explained below, this yields the following pointwise limit of $\munubar(x)$ as $\nu \to 0$:
\begin{equation}
\label{eqn:plim1}
\lim_{\nu \rightarrow 0} \munubar(x) = 0 \quad \text { for } x = 0, \, x=L, \qquad \text{ and } \qquad 
\lim_{\nu \rightarrow 0} \munubar(x) = 1 \quad \text { for } x \in (0,L),
\end{equation}
consistent with the limiting behaviour of minimizers  shown in Figure \ref{fig:ExpCompare}(a). 

The pointwise limit at $x=0$ can be inferred immediately from \eqref{eqn:SSforms} and $c_1 \to 0$, $c_2 \to -1$ as $\nu \to 0$ (Figure \ref{fig:MinConsts}(a)). Furthermore, since $c_2$ approaches a finite value as $\nu \to 0$, at strictly positive $x$ in the support of $\munubar$ we have:
\begin{equation}
\label{eqn:plim2}
\lim_{\nu \rightarrow 0} \munubar(x) = \lim_{\nu \rightarrow 0} c_1 e^{\frac{x}{\sqrt{2\nu}}} + 1, \qquad \text{ for } x \in (0,L].
\end{equation}

In Figure~\ref{fig:MinConsts}(b) we explore the behaviour of $c_1 e^{\frac{x}{\sqrt{2\nu}}}$ as $\nu$ tends to zero, at $x=L$ and $x=L-10^{-2}$, that is, at the end of the support, as well as very close to it. We find the two pointwise limits to be $-1$ and $0$ respectively. From this observation and \eqref{eqn:plim2} we conclude \eqref{eqn:plim1}, also noting that if $\lim_{\nu \rightarrow 0} c_1 e^{\frac{x}{\sqrt{2\nu}}} = 0$ for $x$ arbitrarily close to $L$ then the limit is also zero for all $0 < x < L$. 

Finally, it should be remarked that the values of $\nu$ for which we have calculated numerically $c_1$, $c_2$, and $L$ are not coincidental. These are all the values (less than $1$) for which we can reasonably solve the system of nonlinear equations to find the constants. As $\nu$ decreases one finds that the condition number of the system becomes unmanageable beyond $\nu = 10^{-5}$, 
when the numerical method fails and defects in the solution profile are visibly apparent. It should be also noted that while the system becomes ill-conditioned, any results shown in the paper have been compared satisfactorily versus the results from other methods, namely the finite volume method, so that we are confident in what has been reported here.

\bibliographystyle{abbrv}
\def\url#1{}
\bibliography{lit}

\end{document}